\documentclass[final]{siamltex}
\usepackage{amsmath, amssymb, algorithm, algorithmic}
\usepackage{hyperref}
\usepackage{xifthen}

\newcommand{\E}{\ensuremath{\mathbb{E}}}
\newcommand{\R}{\ensuremath{\mathbb{R}}}

\newcommand{\argmin}[1]{\ensuremath{\underset{#1}{\operatorname{argmin}}}~}
\def\ds{\displaystyle}

\def\cQ{{\cal Q}}
\def\cC{{\cal C}}
\def\cB{{\cal B}}
\def\gt{\ensuremath{\hat{\cal G}(x\md)}}
\def\pt{\ensuremath{\hat{\cal K}_x(\bar x_t)}}
\def\rt{\ensuremath{\hat{\cal K}_y(y_{t+1})}}
\def\gi{\ensuremath{\hat{\cal G}(x_i^{md})}}
\def\pi{\ensuremath{\hat{\cal K}_x(\bar x_i)}}
\def\ri{\ensuremath{\hat{\cal K}_y(y_{i+1})}}
\def\tn{\ensuremath{_{t+1}}}

\def\prob{{\rm Prob}}
\def\itemeqn{\item\abovedisplayskip=2pt\abovedisplayshortskip=0pt~\vspace*{-\baselineskip}}
\newcommand{\ag}[1][]{\ensuremath{
	\ifthenelse{\isempty{#1}}
		{_{t}^{ag}}
		{_{t+#1}^{ag}}}}
\newcommand{\md}[1][]{\ensuremath{
		\ifthenelse{\isempty{#1}}
		{_{t}^{md}}
		{_{t+#1}^{md}}}}
\newcommand{\bg}[1][]{\ensuremath{
		\ifthenelse{\isempty{#1}}
		{\beta_{t}\gamma_t}
		{(\beta_{t}-1)\gamma_t}}}
\newcommand{\rst}[1][]{\ensuremath{
		\ifthenelse{\isempty{#1}}
		{\frac 1{2\eta_t}}
		{\frac 1{2\eta_{t+1}}}}}
\newcommand{\rtt}[1][]{\ensuremath{
		\ifthenelse{\isempty{#1}}
		{\frac 1{2\tau_t}}
		{\frac 1{2\tau_{t+1}}}}}
\newcommand{\xdiff}[1][]{\ensuremath{
		\ifthenelse{\isempty{#1}}
		{\frac 1{2\eta_t}\|x-x_t\|^2 - \frac 1{2\eta_t}\|x-x_{t+1}\|^2}
		{\frac {\gamma_t}{2\eta_t}\|x-x_t\|^2 - \frac {\gamma_t}{2\eta_t}\|x-x_{t+1}\|^2}}}
\newcommand{\ydiff}[1][]{\ensuremath{
		\ifthenelse{\isempty{#1}}
		{\frac 1{2\tau_t}\|y-y_t\|^2 - \frac 1{2\tau_t}\|y-y_{t+1}\|^2}
		{\frac {\gamma_t}{2\tau_t}\|y-y_t\|^2 - \frac {\gamma_t}{2\tau_t}\|y-y_{t+1}\|^2}}}

 \newcommand{\bbe}{\mathbb{E}}
\newcommand{\bbr}{\mathbb{R}}
\def\eqnok#1{(\ref{#1})}
\def\cG{{\cal G}}
\def\hcG{{\hat {\cal G}}}

\def\hcK{{ \hat {\cal K}}}
\def\SO{{\cal SO}}
\def\prob{\mathop{\rm Prob}}

\def\vgap{\vspace*{.1in}}
\newcommand{\beq}{\begin{equation}}
\newcommand{\eeq}{\end{equation}}
\newcommand{\beqa}{\begin{eqnarray}}
\newcommand{\eeqa}{\end{eqnarray}}
\newcommand{\beqas}{\begin{eqnarray*}}
\newcommand{\eeqas}{\end{eqnarray*}}

\title{Optimal Primal-Dual Methods for a Class of Saddle Point Problems}
\author{Yunmei Chen\thanks{Department of Mathematics, University of Florida ({\tt yun@math.ufl.edu}). This author was partially supported by
    NSF grant DMS-1115568} \and
Guanghui Lan\thanks{Department of Industrial and System Engineering, University of Florida ({\tt glan@ise.ufl.edu}).
This author was partially supported by NSF grant CMMI-1000347, ONR grant N00014-13-1-0036 and NSF CAREER Award CMMI-1254446.} \and
Yuyuan Ouyang\thanks{Department of Mathematics, University of Florida ({\tt ouyang@ufl.edu})}}

\begin{document}

\maketitle

\begin{abstract}
We present a novel accelerated primal-dual (APD) method for solving a
class of deterministic and stochastic saddle point problems (SPP). The basic idea of
this algorithm is to incorporate a multi-step acceleration scheme into the primal-dual method
without smoothing the objective function. For deterministic SPP, the APD method achieves the
same optimal rate of convergence as Nesterov's  smoothing technique.  Our stochastic APD method
exhibits an optimal rate of convergence for stochastic SPP not only in terms of its dependence on the
 number of  the iteration, but also on a variety of problem parameters. To the best
of our knowledge, this is the first time that such an optimal algorithm has been developed for stochastic SPP
in the literature. Furthermore,  for both deterministic and stochastic SPP,  the developed APD algorithms
can deal with the situation when the feasible region is unbounded, as long as a saddle point exists.
In the unbounded case, we incorporate the modified termination criterion introduced by Monteiro and
Svaiter in solving SPP problem posed as monotone inclusion, and demonstrate that the rate of convergence of
the APD method depends on the distance from the initial point to the set of optimal solutions.
\end{abstract}


\section{Introduction}

Let ${\cal X}$ and ${\cal Y}$ denote the finite-dimensional vector spaces
 equipped with an inner product
$\langle\cdot,\cdot\rangle$ and norm $\|\cdot\|$, and $X \subseteq {\cal X}$, $Y \subseteq {\cal Y}$ be given closed convex sets.
The basic problem of interest in this paper is the saddle-point problem (SPP) given in the form of:
\begin{equation}
	\label{eqnDSPP}
	\min_{x\in X} \left\{ f(x):= \max_{y\in Y}G(x)+ \langle Kx, y \rangle-J(y)\right\}.
\end{equation}
Here, $G(x)$ is a general smooth convex function such that, for some $L_G \ge 0$,
\begin{align}
	\label{eqnLG}
G(y) - G(x) - \langle \nabla G(x), y - x \rangle \le \frac{L_G}{2} \|y - x\|^2, \ \ \forall x, y \in X,
\end{align}
$K: {\cal X}\rightarrow {\cal Y}$ is  a  linear operator with induced norm
$L_K=\|K\|$, and $J: Y \rightarrow \R $ is a relatively simple, proper,
convex, lower semi-continuous (l.s.c.) function (i.e., problem \eqnok{eqnACPDyt1} is easy to solve).
In particular, if $J$ is the convex conjugate of some convex function $F$ and $Y \equiv {\cal Y}$, then
\eqnok{eqnDSPP} is equivalent to the primal problem:
\begin{equation}
	\label{pProblem}
	\min_{x\in X}G(x) + F(Kx).
\end{equation}
Problems of these types have recently found many applicaitons in data analysis, especially
in imaging processing and machine learning. In many of these applications, $G(x)$
is a convex data fidelity term, while $F(Kx)$ is a certain regularization,
e.g., total variation~\cite{rudin1992nonlinear}, low rank tensor~\cite{KolBad09-1,TSHK11-1},
overlapped group lasso~\cite{JOV09,MaJeObBa11-1}, and graph regularization~\cite{JOV09,TSRZK05-1}.

This paper focuses on first-order methods for solving both determinisitc SPP, where exact first-order information
on $f$ is available, and stochastic SPP, where we only have access to inexact information about
$f$. Let us start by reviewing a few existing first-order methods in both cases.

\subsection{Deterministic SPP}
\label{secDSPPIntro}
Since the objective function $f$ defined in \eqnok{eqnDSPP} is nonsmooth in general,
traditional nonsmooth optimization methods, e.g., subgradient methods,
would exhibit an ${\cal O}(1/\sqrt{N})$ rate of convergence when applied to \eqnok{eqnDSPP}~\cite{nemirovski1983problem},
where $N$ denotes the number of iterations.
However, following the breakthrough paper by Nesterov~\cite{nesterov2005smooth},
much research effort has been devoted to the development of more efficient methods for solving problem \eqnok{eqnDSPP}.

\vgap

\noindent {(1) \sl Smoothing techniques.} In~\cite{nesterov2005smooth}, Nesterov proposed to approximate
the nonsmooth objective function $f$ in \eqnok{eqnDSPP} by a smooth one with Lipschitz-continuous gradient.
Then, the smooth approximation function is minimized by an accelerated gradient
method in \cite{nesterov1983method, nesterov2004introductory}. Nesterov demonstrated in~\cite{nesterov2005smooth}
that, if $X$ and $Y$ are compact, then the rate of convergence of this smoothing scheme applied to \eqnok{eqnDSPP}
can be bounded by:
\begin{equation}\label{eqnOptRateDSPP}
{\cal O} \left(\frac{L_G}{N^2}+\frac{L_K}{N}\right),
\end{equation}
which significantly improves the previous bound ${\cal O}(1/\sqrt{N})$. It can be seen that the rate of convergence in \eqref{eqnOptRateDSPP}
is actually optimal, based on the following observations:
\begin{enumerate}
	\item [a)] There exists a function $G$ with Lipschitz continuous gradients,
	such that for any first-order method, the rate of convergence for solving $\ds\min_{x\in X}G(x)$ is at most $\ds{\cal O}
	\left(L_G/N^2\right)$ \cite{nesterov2004introductory}.
	\item [b)] There exists $b\in Y$, where $Y$ is a convex compact set of $\bbr^m$ for some $m>0$, and a linear bounded operator $K$,
	such that for any first-order method, the rate of convergence for solving
	$\ds\min_{x\in X}\max_{y\in Y}\langle Kx, y\rangle -J(y):=\min_{x\in X}\max_{y\in Y}\langle Kx-b, y\rangle$
	is at most $\ds {\cal O}\left(L_K/N\right)$ \cite{nemirovsky1992information, nemirovski2004prox}.
\end{enumerate}
Nesterov's smoothing technique has been extensively studied, see,
e.g.,~\cite{nesterov2005excessive,AuTe06-1,LaLuMo11-1,dasp08-1,pena08-1,tseng2008accelerated,BeBoCa09-1,Lan13-1}).
Observe that in order to properly apply these smoothing technqiues, we need to assume either $X$ or $Y$ to be bounded. 

\vgap

\noindent{(2) \sl Primal-dual methods.} While Nesterov's smoothing scheme or its variants rely on
a smooth approximation to the orginal problem \eqnok{eqnDSPP}, primal-dual methods work directly
with the original saddle-point problem. This type of method was first presented by Arrow et al.
\cite{arrow1958studies} and named as the primal-dual hybrid gradient (PDHG)
method in \cite{zhu2008efficient}. The results in \cite{zhu2008efficient,chambolle2011first,esser2010general}
showed that the PDHG algorithm, if employed with well-chosen stepsize policies, exhibit very fast convergence in practice,
especially for some imaging applications. Recently Chambolle and Pork~\cite{chambolle2011first}
presented a unified form of primal-dual algorithms, and demonstrated that,
with a properly specified stepsize policy and averaging scheme, these algorithms can also achieve
the ${\cal O} (1/N)$ rate of convergence. They also discussed possible ways to extend primal-dual algorithms
to deal with the case when either $X$ and $Y$ are unbounded. In the original work of
Chambolle and Pork, they assume $G$ to be relatively simple so that the subproblems
can be solved efficiently. With little additional effort, one can show that, by linearizing $G$
at each step, their method can also be applied for a general smooth convex funtion $G$ and
the rate of convergence of this modified algorithm is given by
\beq \label{nemrate}
{\cal O} \left(\frac{L_G + L_K}{N}\right).
\eeq
It should be noted, however, that although both bounds in \eqnok{eqnOptRateDSPP} and \eqnok{nemrate} are ${\cal O}(1/N)$,
the one in \eqnok{eqnOptRateDSPP} has a significantly better dependence on $L_G$.
More specifically, Nesterov's smoothing scheme would allow a very large Lipschitz
constant $L_G$ (as big as ${\cal O}(N)$) without affecting the rate of convergence (up to a constant factor of $2$).
This is desirable in many data analysis applications (e.g., image processing), where
$L_G$ is usually significantly bigger than $L_K$.
Note that the primal-dual methods are also related to
the Douglas-Rachford splitting method \cite{douglas1956numerical} and a pre-conditioned version of
the alternating direction method of multipliers \cite{gabay1976dual}.

\vgap

\noindent{ (3) \sl Extragradient methods for variation inequality (VI) reformulation.}
Motivated by Nesterov's work, Nemirovski presented a mirror-prox method, by modifying
Korpelevich's extragradient algorithm~\cite{Korpelevich1983extrapolation}, for solving a more general class of
variational inequalities \cite{nemirovski2004prox} (see also \cite{JudNem11-2}).
Similar to the primal-dual methods mentioned above, the extragradient methods update iterates on
both the primal space ${\cal X}$ and dual space ${\cal Y}$, and do not require any smoothing technique.
The difference is that each iteration of the extragradient methods requires an extra gradient descent step.
Nemirovski's method, when specialized to \eqnok{eqnDSPP}, also
exhibits a rate of convergence given by \eqnok{nemrate}, which, in view of our
previous discussion, is not optimal in terms of its dependence on $L_G$. It can be shown that,
in some special cases (e.g., $G$ is quadratic),
one can write explicitly the (strongly concave) dual function of $G(x)$ and obtain a result
similar to \eqnok{eqnOptRateDSPP}, e.g., by applying an improved algorithm in \cite{JudNem11-2}.
However, this approach would increase the dimension of the problem and
cannot be applied for a general smooth function $G$. It should be noted that,
while Nemirovski's initial work only considers the case when both $X$ and $Y$ are bounded,
Monteiro and Svaiter~\cite{MonSva09-1} recently showed that extragradient methods
can deal with unbounded sets $X$ and $Y$ by using a slightly modified termination criterion.

\subsection{Stochastic SPP}
While determinisitc SPP has been extensively explored, the study
on stochastic first-order methods for stochastic SPP is still quite limited.
In the stochastic setting, we assume that there exists a {\sl stochastic
oracle} ($\SO$) that can provide unbiased estimators to the gradient operators $\nabla G(x)$ and $(-K x, K^T y)$.
More specifically, at the $i$-th call to $\SO$, $(x_i, y_i) \in X \times Y$ being the input,
the oracle will output the {\sl stochastic gradient}
$(\hcG(x_i), \hcK_x(x_i), \hcK_y(y_i)) \equiv
({\cal G}(x_i, \xi_i), {\cal K}_x(x_i, \xi_i), {\cal K}_y(y_i,\xi_i))$ such that
\beq
\bbe[\hcG(x_i)] = \nabla G(x_i),
 \ \ \ \ \
\bbe\left[ \left(
\begin{array}{c}
-\hcK_x(x_i)\\
\hcK_y(y_i)
\end{array}\right)
\right] = \left(
\begin{array}{c}
-K x_i\\
K^T y_i
\end{array}
\right).
\eeq
Here $\{\xi_i \in \bbr^d\}_{i=1}^\infty$ is a sequence of i.i.d. random variables.
In addition, we assume that, for some $\sigma_{x, G}, \sigma_y, \sigma_{x, K} \ge 0$, the following assumption holds:

\begin{enumerate}
	\renewcommand{\theenumi}{{\bf A\arabic{enumi}}}
	\item 
	\label{itmVar}
	$
	\bbe[\|\hcG(x_i) - \nabla G(x_i)\|_*^2] \le \sigma_{x, G}^2,\ 
	$
	$
	\bbe[\|\hcK_x(x_i) - K x_i\|_*^2] \le \sigma_y^2
	\ \ \
	\mbox{and} \ \ \
	\bbe[\|\hcK_y(y_i) - K^T y_i\|_*^2] \le \sigma_{x, K}^2.
	$
	\newcounter{enuAssumptions}
	\setcounter{enuAssumptions}{\value{enumi}}	
\end{enumerate}

Sometimes we simply denote $\sigma_x := \sqrt{\sigma_{x, G}^2 + \sigma_{x, K}^2}$ for the sake of notational convenience.
Stochastic SPP often appears in machine learning applications. For example,
for problems given in the form of \eqnok{pProblem}, $G(x)$ (resp. $F(Kx)$) can be used to
denote a smooth (resp. nonsmooth) expected convex loss function.
It should also be noted that deterministic SPP is a special case of the above setting
with $\sigma_x = \sigma_y = 0$.

In view of the classic complexity theory for convex programming~\cite{nemirovski1983problem,juditsky2008solving},
a lower bound on the rate of convergence for solving stochastic SPP is given by
\beq \label{lower_bnd}
\Omega\left(\frac{L_G}{N^2} + \frac{L_K}{N} + \frac{\sigma_x + \sigma_y}{\sqrt{N}}\right),
\eeq
where the first two terms follow from the discussion after \eqnok{eqnOptRateDSPP} and
the last term follows from Section 5.3 and 6.3 of \cite{nemirovski1983problem}.
However, to the best of our knowledge, there does not exist an optimal algorithm in the literature which 
exhibits exactly the same rate of convergence as in \eqnok{lower_bnd}, although 
there are a few general-purpose stochastic optimization algorithms which possess different nearly optimal
rates of convergence when applied to above stochastic SPP.

\vgap

\noindent {(1) \sl Mirror-descent stochastic approximation (MD-SA).} The MD-SA method developed
by Nemirovski et al. in \cite{nemirovski2009robust} originates from the classical stochastic approximation (SA) of
Robbins and Monro~\cite{RobMon51-1}. The classical SA mimics the simple gradient descent method by replacing exact gradients
with stochastic gradients, but can only be applied to solve strongly convex problems (see also Polyak \cite{pol90} and Polyak and Juditsky
\cite{pol92}, and Nemirovski et al. \cite{nemirovski2009robust} for an account for the earlier
development of SA methods). By properly modifying the
classical SA, Nemirovski et al. showed in \cite{nemirovski2009robust} that the MD-SA method
can optimally solve general nonsmooth stochastic programming problems. The rate
of convergence of this algorithm, when applied to the stochastic SPP,
is given by (see Section 3 of \cite{nemirovski2009robust})
\[
{\cal O} \left\{(L_G + L_K + \sigma_x + \sigma_y)\frac{1}{\sqrt{N}} \right\}.
\]
However, the above bound is significantly worse than the lower bound in \eqnok{lower_bnd}
in terms of its dependence on both $L_G$ and $L_K$.

\vgap

\noindent {(2) \sl Stochastic mirror-prox (SMP).} In order to improve the convergence
of the MD-SA method, Juditsky et al.~\cite{juditsky2008solving} developed a stochastic
counterpart of Nemirovski's mirror-prox method for solving general variational inequalities.
The stochastic mirror-prox method, when specialized to the above stochastic SPP, yields a rate of convergence given by
\[
{\cal O} \left\{\frac{L_G + L_K}{N} + \frac{\sigma_x + \sigma_y}{\sqrt{N}} \right\}.
\]
Note however, that the above bound is still significantly worse than the lower bound in \eqnok{lower_bnd}
in terms of its dependence on $L_G$.

\vgap

\noindent {(3) \sl Accelerated stochastic approximation (AC-SA).} More
recently, Lan presented in~\cite{lan2012optimal} (see also \cite{GhaLan10-1a,GhaLan10-1b}) a unified optimal method for solving
smooth, nonsmooth and stochastic optimization by developing
a stochastic verstion of Nesterov's method~\cite{nesterov1983method,nesterov2004introductory}.
The developed AC-SA algorithm in~\cite{lan2012optimal}, when applied to the aforementioned
stochastic SPP, possesses the rate of convergence given by
\[
{\cal O} \left\{ \frac{L_G}{N^2} + (L_K + \sigma_x + \sigma_y)\frac{1}{\sqrt{N}}
\right\}.
\]
However, since the nonsmooth term in $f$ of \eqnok{eqnDSPP} has certain special structure,
the above bound is still significantly worse than the lower bound in \eqnok{lower_bnd}
in terms of its dependence on $L_K$. It should be noted that some improvement for AC-SA
has been made by Lin et al. \cite{LinChePen11-1} by applying the smoothing technique
to \eqnok{eqnDSPP}. However, such an improvement works only for the case when $Y$ is bounded
and $\sigma_y = \sigma_{x, K} = 0$. Otherwise, the rate of convergence
of the AC-SA algorithm will depend on the ``variance'' of the stochastic 
gradients computed for the smooth approximation problem, which is usually unknown
and difficult to characterize (see Section 3 for more discussions). 

\vgap

Therefore, none of the stochastic optimization algorithms mentioned above could achieve the lower bound
on the rate of convergence in \eqnok{lower_bnd}. 

\subsection{Contribution of this paper}
Our contribution in this paper mainly consists of the following three aspects.
Firstly, we present a new primal-dual type method, namely the accelerated primal-dual (APD) method, that can achieve the optimal
rate of convergence in \eqnok{eqnOptRateDSPP} for deterministic SPP.
The basic idea of this algorithm is to incorporate a multi-step acceleration scheme into the primal-dual method in \cite{chambolle2011first}.
We demonstrate that, without requiring the application of the smoothing technique, this method can also achieve
the same optimal rate of convergence as Nesterov's smoothing scheme when applied to \eqnok{eqnDSPP}.
We also show that the cost per iteration for APD is comparable to that of Nesterov's smoothing scheme.
Hence our method can efficiently solve problems with a big Lipschtiz constant $L_G$.

Secondly, in order to solve stochastic SPP, we develop a stochastic counterpart of the APD method, namely stochastic APD and
demonstrate that it can actually achieve the lower bound on the rate of convergence in \eqnok{lower_bnd}.
Therefore, this algorithm exhibits an optimal rate of convergence for stochastic SPP not only
in terms of its dependence on $N$, but also on a varity of problem parameters including, $L_G$, $L_K$, $\sigma_x$ and
$\sigma_y$. To the best of our knowledge, this is the first time that such an optimal
algorithm has been developed for stochastic SPP in the literature.
In addition, we investigate the stochastic APD method in more details, e.g., by developing
the large-deviation results associated with the rate of convergence of the stochastic APD method.

Finally, for both deterministic and stochastic SPP, we demonstrate that the developed
APD algorithms can deal with the situation when either $X$ or $Y$ is unbounded,
as long as a saddle point of problem \eqnok{eqnDSPP} exists. We incorporate into the APD method 
the termination criterion employed by Monteiro and Svaiter~\cite{monteiro2011complexity} for solving variational inequalities,
and generalize it for solving stochastic SPP. In both deterministic and stochastic cases, the rate of convergence of the APD 
algorithms will depend on the distance from the initial point to the set of optimal solutions.
\subsection{Organization of the paper} We present the APD methods and discuss their main convergence
properties for solving deterministic and stochastic
SPP problems, respectively, in Sections 2 and 3. 
In order to facilitate the readers, we put the proofs of our main results in 
Section 4. Some brief concluding remarks are made in Section 5.

\section{Accelerated Primal-Dual Methods for Deterministic SPP}
\label{secACPD}
Our goal in this section is to present an accelerated primal-dual method for deterministic SPP and 
discuss its main convergence properties. 

The study on first-order primal-dual method for nonsmooth convex optimization has
 been mainly motivated by solving total variation based image processing problems
(e.g. \cite{zhu2008efficient, esser2010general, pock2009algorithm, chambolle2011first, bonettini2012convergence, he2012convergence}). 
Algorithm~\ref{algCP} shows a primal-dual method summarized in~\cite{chambolle2011first} for solving a special
case of problem (1.1), where $Y=R^m$ for some $m>0$, and $J(y)= F^*(y)$ is the convex conjugate of a convex
and l.s.c. function $F$.
%

\begin{algorithm}
	\caption{Primal-dual method for solving deterministic SPP}
	\label{algCP}
	\begin{algorithmic}[1]
\STATE Choose $x_1\in X$, $y_1\in Y$. Set $\bar{x}_1=x_1$.

\STATE For $t=1, \ldots, N$, calculate 
		\begin{align}
\label{eqnCPy} y_{t+1} & = \argmin{y\in Y}\langle -K\bar{x}_{t}, y \rangle + J(y) + \frac{1}{2\eta_t}\|y-y_t\|^2,\\
\label{eqnCPx}			x_{t+1} & = \argmin{x\in X}G(x) + \langle Kx,  y_{t+1}\rangle + \frac {1}{2\tau_t}\|x-x_t\|^2,\\
		\bar x_{t+1} & = \theta_t(x_{t+1} - x_{t}) + x_{t+1}.
		\end{align}
		\STATE Output $\ds x^N=\frac 1N\sum_{t=1}^N x_t$, $\ds y^N=\frac 1N\sum_{t=1}^N y_t$.
	\end{algorithmic}
\end{algorithm}

\begin{algorithm}
	\caption{Accelerated primal-dual method for deterministic SPP}
	\label{algACPD}
	\begin{algorithmic}[1]
		\STATE Choose $x_1\in X, y_1\in Y$. Set $x_1^{ag}=x_1, y_1^{ag}=y_1$, $\bar x_1 = x_1$.
		\STATE For $t=1,2,\ldots, N-1$, calculate
		\begin{align}
			\label{eqnACPDxmd}x\md & = (1-\beta_t^{-1})x\ag + \beta_t^{-1}x_t,\\
			\label{eqnACPDyt1}y_{t+1} & = \argmin{y\in Y}\langle -K\bar x_t,y\rangle + J(y) + \frac 1{\tau_t}V_Y(y, y_t),\\
			\label{eqnACPDxt1}x_{t+1} & = \argmin{x\in X}\langle \nabla G(x\md), x\rangle + \langle x, K^Ty\tn\rangle + \frac 1{\eta_t}V_X(x, x_t),\\			
			\label{eqnACPDxag}x\ag[1] & = (1-\beta_t^{-1})x\ag + \beta_t^{-1}x_{t+1},\\
			\label{eqnACPDyag}y\ag[1] & = (1-\beta_t^{-1})y\ag + \beta_t^{-1}y_{t+1},\\
			\label{eqnACPDxbar}		\bar x_{t+1} & = \theta_{t+1}(x_{t+1} - x_{t}) + x_{t+1}.			
		\end{align}
		\STATE Output $x_{N}^{ag},y_{N}^{ag}$.
	\end{algorithmic}
\end{algorithm}

The convergence of the sequence $\{(x_t, y_t)\}$ in Algorithm \ref{algCP} has been studied in \cite{pock2009algorithm, esser2010general, chambolle2011first,bonettini2012convergence, he2012convergence} for various choices of $\theta_t$, and under different conditions on
the  stepsizes  $\tau_t$ and $ \eta_t $. 
However, the rate of convergence for this algorithm has only been discussed by Chambolle and Pock in \cite{chambolle2011first}. 
More specifically, they assume that the constant stepsizes are used, i.e., $\tau_t = \tau$, $\eta_t = \eta$
and $\theta_t = \theta$ for some $\tau, \eta, \theta > 0$ for all $t \ge 1$. 
If $\tau \eta L_K^2 <1$, where $L_K=\|K\|$, then the output $(x^N, y^N)$ possesses 
a rate of convergence of ${\cal O}(1/N)$ for $\theta=1$, and of $O(1/\sqrt{N})$ for $\theta=0$, in
terms of partial duality gap (duality gap in a bounded domain, see \eqnok{def_gap} below).

One possible limitation of \cite{chambolle2011first} is that both $G$ and $J$ need to be simple enough 
so that the two subproblems \eqnok{eqnCPy} and \eqnok{eqnCPx} in Algorithm~\ref{algCP} are easy to solve.
 To make Algorithm \ref{algCP} applicable to more practical problems we  consider more general cases, where $J$ is simple, but $G$ may not be so.  
 In particular, we assume that $G$ is a general smooth convex function satisfying \eqnok{eqnDSPP}.
 In this case, we can replace $G$ in \eqref{eqnCPx} by its linear approximation $G(x_t) + \langle\nabla G(x_t), x - x_t\rangle$.
 Then \eqref{eqnCPx} becomes
\begin{equation} \label{eqnCPxl}
    x_{t+1}  = \argmin{x\in X}\langle\nabla G(x_t), x\rangle + \langle Kx,  y_{t+1}\rangle + \frac {1}{2\tau}\|x-x_t\|^2.
\end{equation}
In the following context, we will refer to this modified algorithm as the ``linearized version'' of Algorithm \ref{algCP}. 
By some extra effort we can show that, if  for $t=1, \ldots, N$,
$0<\theta_t = \tau_{t-1}/\tau_t=\eta_{t-1}/\eta_t \leq 1 $, and
$L_G \tau_t+L_K^2\eta_t\tau_{t}\leq 1$, then  $(x^N,y^N)$ has an ${\cal O}((L_G+L_K)/N)$ rate of convergence in the sense of the partial duality gap.

As discussed in Section 1, the aforementioned rate of convergence for the linearized version of Algorithm \ref{algCP}
is the same as that proved in \cite{chambolle2011first}, and not optimal
in terms of its dependence on $L_G$ (see \eqref{nemrate}). However, this algorithm solves the problem (1.1) directly without smoothing the nonsmooth objective function. Considering the primal-dual method as an alternative to Nestrov's smoothing method, and inspired by his idea of using accelerated gradient descent algorithm to solve the smoothed problem \cite{nesterov1983method, nesterov2004introductory, nesterov2005smooth}, we propose
the following accelerated primal-dual algorithm that integrates the accelerated gradient descent algorithm into
the linearized version of Algorithm \ref{algCP}.

Our accelerated primal-dual (APD) method is presented in Algorithm \ref{algACPD}.
Observe that in this algorithm, the superscript ``ag'' stands for ``aggregated'', and ``md'' stands for ``middle''. 
The functions $V_X(\cdot, \cdot)$ and $V_Y(\cdot, \cdot)$ are Bregman divergences defined as
\begin{align}
	\label{eqnVx}
	V_X(x, u) & :=d_X(x)-d_X(u)-\langle\nabla d_X(u), x-u\rangle, \forall x, u\in X,\\
	\label{eqnVy}
	V_Y(y, v) & :=d_Y(y)-d_Y(v)-\langle\nabla d_Y(v), y-v\rangle, \forall y, v\in Y,
\end{align}
where $d_X(\cdot)$ and $d_Y(\cdot)$ are strongly convex functions with strong convexity parameters $\alpha_X$ and $\alpha_Y$. For example, under the Euclidean setting, we can simply set $V_X(x, x_t): = \|x - x_t\|^2/2$ and $V_Y(y, y_t):= \|y - y_t\|^2/2$, and $\alpha_X=\alpha_Y=1$. We assume that $J(y)$ is a simple convex function, so that the optimization problem in \eqref{eqnACPDyt1} can be solved efficiently.

Note that if $\beta_t = 1$ for all $t \ge 1$, then $x_{t}^{md}=x_t$, $x_{t+1}^{ag}=x_{t+1}$,
and  Algorithm \ref{algACPD} is the same as the linearized version of Algorithm \ref{algCP}.
However, by specifying a different selection of $\beta_t$ (e.g., $\beta_t = O(t)$), 
we can significantly improve the rate of convergence of Algorithm~\ref{algACPD} in terms
of its dependence on $L_G$. It should be noted that the iteration cost for the APD algorithm
is about the same as that for  the linearized version of Algorithm \ref{algCP}.

In order to analyze the convergence of Algorithm~\ref{algACPD}, it is necessary to 
introduce a notion to characterize the solutions of \eqref{eqnDSPP}. Specifically,
denoting $Z = X\times Y$, for any $\tilde z = (\tilde x, \tilde y) \in Z$ and $z = (x, y)\in Z$, 
we define
	\begin{equation}
		\label{eqnQ}
		Q(\tilde z, z):=\left[G(\tilde x)+\langle K\tilde x, y\rangle-J(y)\right] - \left[G(x)+\langle Kx, \tilde y\rangle-J(\tilde y)\right].
	\end{equation}
It can be easily seen that $\tilde z$ is a solution of problem \eqref{eqnDSPP},
if and only if $Q(\tilde z, z)\leq 0$ for all $z\in Z$.
Therefore, if $Z$ is bounded, it is suggestive to use the gap function
\beq \label{def_gap}
g(\tilde z) := \max_{z \in Z} Q(\tilde z, z)
\eeq
to assess the quality of a feasible solution $\tilde z \in Z$. In fact, we can show that
$f(\tilde x) - f^*\leq g(\tilde z)$ for all $\tilde z\in Z$,
where $f^*$ denotes the optimal value of problem \eqref{eqnDSPP}. 
However, if $Z$ is unbounded, then $g(\tilde z)$ is not well-defined even for a nearly optimal solution $\tilde z \in Z$.
Hence, in the sequel, we will consider the bounded and unbounded case separately,
by employing a slightly different error measure for the latter situation.

\vgap

The following theorem describes the convergence properties of Algorithm \ref{algACPD} when $Z$ is bounded. 

\begin{theorem}
	\label{thmConvergenceBD}
	Suppose that for some $\Omega_X, \Omega_Y >0$,
	\beq \label{boundness}
	\ds\sup_{x_1, x_2\in X}V_X(x_1, x_2)\leq \Omega_X^2 \ \  \mbox{and} \ \ \ds\sup_{y_1, y_2\in Y}V_Y(x_1, x_2)\leq \Omega_Y^2.
	\eeq 
	Also assume that the 
	parameters $\beta_t, \theta_t, \eta_t, \tau_t$ in Algorithm \ref{algACPD} are chosen such that for all $t\ge 1$,
		\begin{align}
			\label{eqnIneqBetaGamma}
			&\beta_1=1,\ 
			\beta_{t+1}-1 = \beta_t\theta_{t+1},
			\\
			\label{eqnIncStep}
			& 0< \theta_t \leq \min\{\ds\frac{\eta_{t-1}}{\eta_{t}},  \ds\frac{\tau_{t-1}}{\tau_{t}} \},
			\\
			\label{eqnIneqLGLK}&\frac{\alpha_X}{\eta_t}-\frac{L_G}{\beta_t}-\frac{L_K^2\tau_{t}}{\alpha_Y}\geq 0.
		\end{align}
		Then for all $t\ge 1$, 
			\begin{equation}
				\label{eqnErrBDD}
				g(z_{t+1}^{ag})\leq \frac {1}{\beta_t\eta_t}\Omega_X^2 
							+\frac{1}{\beta_t\tau_t}\Omega_Y^2.
			\end{equation}			
\end{theorem}
		
There are various options for choosing the parameters $\beta_t, \eta_t, \tau_t$ and $\theta_t$ 
such that \eqnok{eqnIneqBetaGamma}--\eqnok{eqnIneqLGLK} hold.
Below we provide such an example.

\begin{corollary}
	\label{corStepsizeBD}
	Suppose that \eqnok{boundness} holds. In Algorithm \ref{algACPD}, if the parameters are set to
	\begin{align}
		\label{eqnStepsizeBD}
		\beta_t  =\frac{t+1}{2},\
		\theta_t  = \frac{t-1}{t},\
		 \eta_t  = \frac{\alpha_Xt}{2L_G+tL_KD_Y/D_X} \ \mbox{and} \ \tau_t=\frac{\alpha_YD_Y}{L_KD_X},
	\end{align}	
	where $D_X := \Omega_X \ds\sqrt{2/\alpha_X}$ and $D_Y := \Omega_Y \ds\sqrt{2/\alpha_Y}$,
	then for all $t \ge 2$,
		\begin{align}
			\label{eqnRateBD}
			&\ g(z_{t}^{ag}) \leq \frac{2L_GD_X^2}{t(t-1)} + \frac{2L_K D_X D_Y}{t}.
		\end{align}	
	\begin{proof}
	It suffices to verify that the parameters in \eqref{eqnStepsizeBD} satisfies \eqref{eqnIneqBetaGamma}--\eqref{eqnIneqLGLK} in Theorem \ref{thmConvergenceBD}. It is easy to check that \eqref{eqnIneqBetaGamma} and \eqref{eqnIncStep} hold. Furthermore,
	\begin{align*}
		\frac{\alpha_X}{\eta_t} - \frac{L_G}{\beta_t} - \frac{L_K^2\tau_t}{\alpha_Y} = \frac{2L_G + tL_KD_Y/D_X}{t} - \frac{2L_G}{t+1} - \frac{L_KD_Y}{D_X} \geq 0,
	\end{align*}
	so \eqref{eqnIneqLGLK} holds. Therefore, by \eqref{eqnErrBDD}, for all $t \ge 1$ we have
	\begin{align*}
		g(z_{t}^{ag})&\leq \frac {1}{\beta_{t-1}\eta_{t-1}}\Omega_X^2 
		+\frac{1}{\beta_{t-1}\tau_{t-1}}\Omega_Y^2 
		= \frac{4L_G + 2(t-1)L_KD_Y/D_X}{\alpha_X t(t-1)}\cdot\frac{\alpha_X}{2}D_X^2 + \frac{2L_KD_X/D_Y}{\alpha_Y t}\cdot\frac{\alpha_Y}{2} D_Y^2
		\\
		& = \frac{2L_GD_X^2}{t(t-1)} + \frac{2L_K D_X D_Y}{t}.
	\end{align*}

	\end{proof}
\end{corollary}

Clearly, in view of \eqnok{eqnOptRateDSPP}, the rate of convergence of Algorithm \ref{algACPD} 
applied to problem \eqnok{eqnDSPP} is optimal when the parameters are chosen 
according to \eqnok{eqnStepsizeBD}. 
Also observe that 
we need to estimate $D_Y/D_X$ to use these parameters. However, 
it should be pointed out that replacing the ratio $D_Y/D_X$ in \eqnok{eqnStepsizeBD} by any positive constant 
only results an increase in the RHS of  \eqnok{eqnRateBD} by a constant factor. 


\vgap

Now, we study the convergence properties of the APD algorithm for the case when $Z=X \times Y$ is unbounded,
by using a perterbation-based termination criterion recently employed by Monteiro and Svaiter and applied to SPP~\cite{MonSva09-1, monteiro2012iteration, monteiro2011complexity}. 
This termination criterion is based on the enlargement of a maximal monotone operator, which is first 
introduced in \cite{burachik1997enlargement}. One advantage of using this criterion is that its definition does not 
depend on the boundedness of the domain of the operator. More specifically, 
as shown in \cite{monteiro2011complexity, MonSva09-1}, there always exists a purterbation vector $v$ such that 
\beq \label{def_gapt}
\tilde g(\tilde z, v) := \max_{z \in Z} Q(\tilde z, z) - \langle v, \tilde z - z \rangle
\eeq
is well-defined, although the value of $g(\tilde z)$ in \eqnok{def_gap} may be unbounded if $Z$ is unbounded. 
In the following result, we show that the APD algorithm can compute a nearly optimal solution $\tilde z$
with a small residue $\tilde g(\tilde z, v)$, for a small purterbation vector $v$ (i.e., $\|v\|$ is small). In addition,
our derived iteration complexity bounds are proportional to the distance from the initial point to the solution set.


\begin{theorem}
	\label{thmConvergenceUB}
	Let $\{z\ag\} = \{(x\ag, y\ag)\}$ be the iterates generated by Algorithm \ref{algACPD} with $V_X(x, x_t) = \|x - x_t\|^2/2$ 
	and $V_Y(y, y_t) = \|y - y_t\|^2/2$. Assume that the parameters $\beta_t, \theta_t, \eta_t$ and $\tau_t$ 
	satisfy \eqref{eqnIneqBetaGamma}, 
	\begin{align}
		\label{eqnFlatStep}
		&\theta_t = \ds\frac{\eta_{t-1}}{\eta_{t}} = \ds\frac{\tau_{t-1}}{\tau_{t}}, 
		\\
		\label{eqnIneqLGLKp}&\frac{\alpha_X}{\eta_t}-\frac{L_G}{\beta_t}-\frac{L_K^2\tau_{t}}{p\alpha_Y}\geq 0,
	\end{align}
for all $t\ge 1$ and for some $0<p<1$,	then there exists a perturbation vector $v_{t+1}$ such that 
				\begin{equation}
					\label{eqnErrUBD}
				\tilde g(z_{\tn}^{ag}, v_{t+1})	\leq \frac{(2-p)D^2}{\beta_t\eta_t(1-p)} =: \varepsilon\tn 
				\end{equation}
for any $t\ge 1$. Moreover, we have
				\begin{align}
					\label{eqnvbound}
					& \|v\tn\|\leq \frac 1{\beta_t\eta_t}\|\hat x - x_1\| + \frac 1{\beta_t\tau_t}\|\hat y - y_1\| + \left[\frac{1}{\beta_t\eta_t}\left(1 + \sqrt{\frac{\eta_1}{\tau_1(1-p)}}\right) + \frac{2L_K}{\beta_t} \right]D,
				\end{align}
where $(\hat x, \hat y)$ is a pair of solutions for problem \eqnok{eqnDSPP}
and 
	\begin{equation}
		\label{eqnD}
		D := \sqrt{\|\hat x - x_1\|^2 + \frac{\eta_1}{\tau_1}\|\hat y - y_1\|^2}.
	\end{equation}
\end{theorem}

Below we suggest a specific parameter setting which satisfies
 \eqnok{eqnIneqBetaGamma}, \eqnok{eqnFlatStep} and \eqnok{eqnIneqLGLKp}.
\begin{corollary}
	\label{corStepsizeUB}
	In Algorithm \ref{algACPD}, if $N$ is given and the parameters are set to 	\begin{align}
		\label{eqnStepsizeUB}
		\beta_t  =\frac{t+1}{2},\
		\theta_t  = \frac{t-1}{t},\ 
		\eta_t  = \frac{t+1}{2(L_G+NL_K)}, \mbox{and} \ \tau_t = \frac{t+1}{2NL_K} 
	\end{align}	
	then there exists $v_N$ that satisfies \eqref{eqnErrUBD} with 
\beq \label{eqnErrUB}
			\varepsilon_N \leq \frac{10L_G\hat D^2}{N^2} + \frac{10L_K\hat D^2}{N} \ \ \mbox{and} \ \
			 \|v_N\| \leq \frac{15L_G\hat D}{N^2} + \frac{16L_K\hat D}{N},			 
\eeq	
	where $\hat D = \sqrt{\|\hat x - x_1\|^2 + \|\hat y - y_1\|^2}$.

	\begin{proof}
		For the parameters $\beta_t$, $\gamma_t$, $\eta_t$, $\tau_t$ in \eqref{eqnStepsizeUB}, it is clear that \eqref{eqnIneqBetaGamma}, \eqref{eqnFlatStep} holds. Furthermore, let $p=1/4$, for any $t=1,\ldots,N-1$, we have
			\begin{align*}
			& \frac{1}{\eta_t}-\frac{L_G}{\beta_t}  - \frac{L_K^2\tau_t}{p} = \frac{2L_G + 2L_KN}{t+1} - \frac{2L_G}{t+1} - \frac{2L_K^2(t+1)}{L_KN} \geq \frac{2L_KN}{t+1} - \frac{2L_K(t+1)}{N}\geq 0,
		\end{align*}
		thus \eqref{eqnIneqLGLKp} holds. By Theorem \ref{thmConvergenceUB}, inequalities \eqref{eqnErrUBD} and \eqref{eqnvbound} hold. Noting that $\eta_t\leq \tau_t$, in \eqref{eqnErrUBD} and \eqref{eqnvbound} we have $D\leq \hat D$, $\|\hat x - x_1\| + \|\hat y - y_1\| \leq \sqrt{2}\hat D$, hence
		\begin{align*}
			\|v\tn\| 
			&\leq  \frac{\sqrt 2\hat D}{\beta_t\eta_t} + \frac{(1 + \sqrt{4/3})\hat D}{\beta_t\eta_t} + \frac{2L_K\hat D}{\beta_t},
			\\ 
			\varepsilon\tn & \leq\frac{(2-p) \hat D^2}{\beta_t\eta_t(1-p)} = \frac{7\hat D^2}{3\beta_t\eta_t}.
		\end{align*}
Also	note that by  \eqnok{eqnStepsizeUB},
		\[\frac 1{\beta_{N-1}{\eta_{N-1}}} = \frac{4(L_G + L_KN)}{N^2} = \frac{4L_G}{N^2} + \frac{4L_K}{N}.\]
Using the above three relations and the definition of $\beta_t$ in \eqnok{eqnStepsizeUB},
we obtain \eqnok{eqnErrUB} after simpliying the constants.
	\end{proof}
		
\end{corollary}

It is interesting to notice that, if the parameters in Algorithm~\ref{algACPD} are set to \eqnok{eqnStepsizeUB}, then
both residues $\varepsilon_N$ and $\|v_N\|$ in \eqnok{eqnErrUB} reduce to zero with approximately the
same rate of convergence (up to a factor of $\hat D$).
Also observe that in Theorem~\ref{thmConvergenceUB} and Corollary~\ref{corStepsizeUB}, we
fix $V_X(\cdot,\cdot)$ and $V_Y(\cdot,\cdot)$ to be regular distance functions rather than more general Bregman divergences.
This is due to fact that we need to apply the Triagular inequality associated with $\sqrt{V_X(\cdot,\cdot)}$
and $\sqrt{V_Y(\cdot,\cdot)}$, while such an inequality does not necessarily hold for Bregman divergences in general.

\section{Stochastic APD Methods for Stochastic SPP}
Our goal in this section is to present a stochastic APD method for stochastic SPP
(i.e., problem \eqnok{eqnDSPP} with a stochastic oracle) and demonstrate that
it can actually achieve the lower bound in \eqnok{lower_bnd} on the rate of convergence for stochastic SPP.

The stochastic APD method is a stochastic counterpart of the APD algorithm in Section 2,
obtained by simply replacing the gradient operators
$- K \bar x_t$, $\nabla G(x_t^{md})$ and $K^T y_{t+1}$, used in \eqnok{eqnACPDyt1} and
\eqnok{eqnACPDxt1}, with the stochastic gradient operators
computed by the $\SO$,  i.e., $-\pt$, $\gt$ an $\rt$, respectively.
This algorithm is formally described as in Algorithm~\ref{algACPDS}.


\begin{algorithm}[!h]
	\caption{Stochastic APD method for stochastic SPP}
	\label{algACPDS}
	\begin{algorithmic}
		\STATE Modify \eqref{eqnACPDyt1} and \eqref{eqnACPDxt1} in Algorithm \ref{algACPD} to
		\begin{align}
			\label{eqnACPDSyt1}y_{t+1} & = \argmin{y\in Y}\langle -\pt,y\rangle + J(y) + \frac 1{\tau_t}V_Y(y, y_t)\\
			\label{eqnACPDSxt1}x_{t+1} & = \argmin{x\in X}\langle \gt, x\rangle + \langle x, \rt\rangle + \frac 1{\eta_t}V_X(x, x_t)
		\end{align}
	\end{algorithmic}
\end{algorithm}

A few more remarks about the development of the above stochastic APD method are in order. 
Firstly, observe that, although primal-dual methods have been extensively studied for
solving deterministic saddle-point problems, it seems that these types of methods
have not yet been generalized for stochastic SPP in the literature. 
Secondly, as noted in Section~1, one possible way to solve stochastic SPP 
is to apply the AC-SA algorithm in \cite{lan2012optimal}
to a certain smooth approximation of \eqnok{eqnDSPP} by Nesterov \cite{nesterov2005smooth}. 
However, the rate of convergence of this approach will depend on the variance of the stochastic gradients
computed for the smooth approximation problem, which is usually unkown and difficult 
to characterize. On the other hand, the stochastic APD method described above
works directly with the original problem without requring the application of the smoothing technique,
and its rate of convergence will depend on the variance of the stochastic gradient operators computed
for the original problem, i.e., $\sigma_{x,G}^2$, $\sigma_y^2$ and $\sigma_{x,K}^2$ in A1.
We will show that it can achieve exactly the lower bound in \eqnok{lower_bnd} on the rate of convergence for stochastic SPP. 

Similarly to Section \ref{secACPD}, we use the two gap functions $g(\cdot)$ and $\tilde g(\cdot, \cdot)$, respectively,
defined in \eqnok{def_gap} and \eqnok{def_gapt} as the termination criteria for the stochastic 
APD algorithm, depending on whether the feasible set $Z = X \times Y$ is bounded or not.
Since the algorithm is stochastic in nature, for both cases we establish its
expected rate of convergence in terms of $g(\cdot)$ or $\tilde g(\cdot,\cdot)$,
i.e., the ``average'' rate of convergence over many runs of the algorithm.
In addition, we show that if $Z$ is bounded, then the convergence of the APD algorithm
can be strengthened under the following ``light-tail'' assumption on $\cal SO$.

	\begin{enumerate}
	\setcounter{enumi}{\theenuAssumptions}	
	\renewcommand{\theenumi}{{\bf A\arabic{enumi}}}
	\itemeqn 
	\label{itmLT}
	\begin{align*}
		&\E\left[\exp\{\|\nabla G(x) - \hat\cG(x)\|_*^2/\sigma_{x, G}^2 \} \right]\leq \exp\{1\}, \ \
		\E\left[\exp\{\|K x - \hat{\cal K}_x(x)\|_*^2/\sigma_y^2 \} \right]\leq \exp\{1\} \\
		& \mbox{and} \ \ \E\left[\exp\{\|K^T y - \hat{\cal K}_y(y)\|_*^2/\sigma_{x, K}^2 \} \right]\leq \exp\{1\}.
	\end{align*}
	\end{enumerate}
It is easy to see that A2 implies A1 by Jensen's inequality.

\vgap

Theorem~\ref{thmConvergenceBDS} below summarizes the convergence properties of  Algorithm \ref{algACPDS} when $Z$ is bounded.
Note that the following quanity will be used in the statement of this result and
the convergence analysis of the APD algorithms (see Section 4):
\beq		\label{eqnGamma}
			\gamma_t = \left\{
				\begin{array}{ll}
					1, & t=1,\\
					\theta_t^{-1} \gamma_{t-1}, & t\ge 2.
				\end{array}
				 \right.
 \eeq

\begin{theorem}
	\label{thmConvergenceBDS}
	Suppose that \eqnok{boundness} holds for some $\Omega_X,\Omega_Y > 0$. Also assume that for all $t \ge 1$, the parameters $\beta_t, \theta_t, \eta_t$
	and $\tau_t$ in Algorithm \ref{algACPDS} satisfy \eqref{eqnIneqBetaGamma},  \eqref{eqnIncStep}, and 
		\begin{align}
			\label{eqnIneqLGLKS}&\frac{q\alpha_X}{\eta_t}-\frac{L_G}{\beta_t}-\frac{L_K^2\tau_{t}}{p\alpha_Y}\geq 0
		\end{align}
		 for some $p, q\in(0,1)$. Then,
		\begin{enumerate}
		\renewcommand{\theenumi}{(\alph{enumi})}
		\item 
		\label{itmEQBoundBDS}
		Under assumption \ref{itmVar}, for all $t\ge 1$,
		\begin{equation}
			\label{eqnEQBoundBDS}
			\begin{aligned}
			\E[g(z\ag[1])]\leq \cQ_0(t),
			\end{aligned}
		\end{equation}
		where
		\beq
			\label{eqnQ0}
			\begin{array}{l}
			\cQ_0(t)  := \frac {1}{\beta_t \gamma_t}
			\left\{ \frac{2\gamma_t}{\eta_t} \Omega_X^2  +  \frac{2\gamma_t}{\tau_t} \Omega_Y^2 \right\}
			 \ +  \\
 \frac{1}{2 \beta_t \gamma_t} \sum_{i=1}^{t}\left\{ \frac{(2-q)\eta_i\gamma_i}{(1-q)\alpha_X} \sigma_x^2 + 
\frac{(2-p)\tau_i\gamma_i}{(1-p)\alpha_Y} \sigma_y^2 
			\right\}	.
			\end{array}
		\eeq
		
		\item 
		\label{itmHighProb}
		Under assumption \ref{itmLT}, for all $\lambda>0$ and $t\ge 1$, 
		\begin{align}
			\label{eqnHighProb}	
			 Prob\{g(z\ag[1]) > {\cal Q}_0(t) + \lambda{\cal Q}_1(t) \} &\leq 3\exp\{-\lambda^2/3\} + 3\exp\{-\lambda\},
		\end{align}
		where
		\beq 			\label{eqnQ1}
			\begin{array}{l}
			{\cal Q}_1(t):= \frac 1{{\beta_t\gamma_t}} \left(\frac{\sqrt 2\sigma_{x}\Omega_X}{\sqrt{\alpha_X}}+\frac{\sigma_y\Omega_Y}{{\sqrt{\alpha_Y}}}\right)\sqrt{2\sum_{i=1}^{t}\gamma_i^2 }
			+  \\
 \frac{1}{2 \beta_t \gamma_t} \sum_{i=1}^{t}\left\{ \frac{(2-q)\eta_i\gamma_i}{(1-q)\alpha_X} \sigma_x^2 + 
\frac{(2-p)\tau_i\gamma_i}{(1-p)\alpha_Y} \sigma_y^2 
			\right\}.
			\end{array}
			\eeq
	\end{enumerate}
\end{theorem}	

We provide below a specific choice of the parameters $\beta_t$, $\theta_t$, $\eta_t$ and $\tau_t$
for the stochastic APD method for the case when $Z$ is bounded.

\begin{corollary}
	\label{corStepsizeBDS}
Suppose that \eqnok{boundness} holds and let $D_X$ and $D_Y$ be defined
in Corolloary~\ref{corStepsizeBD}.
	In Algorithm~\ref{algACPDS}, if $N \ge 1$ is given and the parameters are set to
	\begin{align}
		\label{eqnStepsizeBDS}
		\begin{aligned}
		\beta_t & =\frac{t+1}{2},\
		\theta_t  = \frac{t-1}{t},\
		\eta_t = \frac{2\alpha_XD_X t}{6L_GD_X+3L_KD_Y (N-1)+3\sigma_xN\sqrt{N-1}} \ \mbox{and}\\
		\tau_t & = \frac{2\alpha_YD_Y t}{3L_KD_X (N-1)+3\sigma_yN\sqrt{N-1}}.
		\end{aligned}
	\end{align}
Then under Assumption \ref{itmVar}, we have 
	\begin{align}
		\label{eqnEQBoundBDSN}
		E[g(z^{ag}_{N})]\leq \frac{6L_GD_X^2}{N(N-1)} + \frac{6L_KD_XD_Y}{N} + \frac{4(\sigma_xD_X +\sigma_yD_Y) }{\sqrt{N-1}} 
		=: \cC_{0}(N).
	\end{align}
	If in addition, Assumption \ref{itmLT} holds, then for all $\lambda>0$, we have  
	\begin{align}
		\label{eqnHighProbN}
		Prob\{g(z_N^{ag}) > \cC_0(N) + \lambda\cC_1(N) \} \leq 3\exp\{-\lambda^2/3\} + 3\exp\{-\lambda\},
	\end{align}
	where
	\begin{align}
		\label{eqnC1N}
		\cC_1(N) = \frac{3(\sigma_xD_X+\sigma_yD_Y)}{\sqrt{N-1}}.
	\end{align}
 	\begin{proof}
		First we check that the parameters in \eqnok{eqnStepsizeBDS} satisfy the conditions in Theorem \ref{thmConvergenceBDS}. The inequalities \eqref{eqnIneqBetaGamma} and \eqref{eqnIncStep} can be checked easily. Furthermore, for all $t=1,\ldots,N-1$, setting $p = q = 2/3$ we have
			\begin{align*}
				&\ \frac{q\alpha_X}{\eta_t}-\frac{L_G}{\beta_t}-\frac{L_K^2\tau_{t}}{p\alpha_Y}\geq\frac{2L_GD_X+L_KD_Y (N-1)}{D_X t}-\frac{2L_G}{t+1}-\frac{L_K^2D_Y t}{L_KD_X (N-1)}\geq 0,
			\end{align*}
			thus \eqref{eqnIneqLGLKS} hold, and hence Theorem \ref{thmConvergenceBDS} holds.
		
		To show \eqref{eqnEQBoundBDSN} and \eqref{eqnHighProbN}, it suffices to show that $\cC_0(N)\geq \cQ_0(N-1)$ and $\cC_1(N)\geq\cQ_1(N-1)$. Observe that by \eqref{eqnGamma} and \eqref{eqnStepsizeBDS}, we have $\gamma_t = t$. Also, observe that 
			$\ds \sum_{i=1}^{N-1}i^2 \leq (N-1)N^2 / 3$, 
			thus
			\begin{align*}
				& \frac{1}{\gamma_{N-1}}\sum_{i=1}^{N-1}\eta_i\gamma_i \leq \frac{2\alpha_XD_X}{3(N-1)^{3/2}N\sigma_x}\sum_{i=1}^{N-1}i^2
				\leq
				\frac{2\alpha_XD_XN}{9{\sigma_x\sqrt{N-1}}},
				\\
				& \frac{1}{\gamma_{N-1}}\sum_{i=1}^{N-1}\tau_i\gamma_i \leq \frac{2\alpha_YD_Y}{3(N-1)^{3/2}N\sigma_y}\sum_{i=1}^{N-1}i^2
				\leq 
				\frac{2\alpha_YD_Y N}{9\sigma_y\sqrt{N-1}}.
			\end{align*}
			Apply the above bounds to \eqref{eqnEQBoundBDS} and \eqref{eqnQ1}, we get
			\begin{align*}
				&\ \cQ_0(N-1)
				\\
				\leq &\ \frac{2}{N}\left(\frac{6L_GD_X+3L_KD_Y (N-1)+3N\sqrt{N-1}{\sigma_x}}{\alpha_XD_X(N-1)}\cdot \frac{\alpha_X}{2}D_X^2 
				 + \frac{3L_KD_X (N-1)+3N\sqrt{N-1}\sigma_y}{\alpha_YD_Y(N-1)}\cdot \frac{\alpha_Y}{2}D_Y^2\right.
				\\
				&\ \left.+ \frac{2\sigma_x^2}{\alpha_X}\cdot \frac{2\alpha_XD_XN}{9{\sigma_x\sqrt{N-1}}}
				+\frac{2\sigma_y^2}{\alpha_Y}\cdot \frac{2\alpha_YD_Y N}{9\sigma_y\sqrt{N-1}}\right)
				\leq \cC_0(N),
				\\
				&\ \cQ_1(N-1)
				\\
				\leq &\ 
				\frac{2}{N(N-1)}\left(\sigma_{x}D_X+\frac{\sigma_yD_Y}{{\sqrt 2}}\right)\sqrt{\frac{2(N-1)N^2}{3}} 
				+ \frac{4\sigma_x^2}{\alpha_XN}\cdot \frac{2\alpha_XD_XN}{9{\sigma_x}\sqrt{N-1}}
				+\frac{4\sigma_y^2}{\alpha_YN}\cdot \frac{2\alpha_YD_YN}{9\sigma_y\sqrt{N-1}}
				\leq 
				\cC_1(N),
			\end{align*}
	so \eqref{eqnEQBoundBDSN} and \eqref{eqnHighProbN} holds.
	\end{proof}	

\end{corollary}

\vgap

Comparing the rate of convergence established in \eqnok{eqnEQBoundBDSN} 
with the lower bound in \eqnok{lower_bnd}, we can clearly see that
the stochastic APD algorithm is an optimal method for solving the stochastic saddle-point problems.
More specifically, in view of \eqnok{eqnEQBoundBDSN}, this algorithm allows us to have very large 
Lipschitz constants $L_G$ (as big as ${\cal O}(N^\frac{3}{2})$) and $L_K$
(as big as ${\cal O}(\sqrt{N})$) without significantly affecting its rate of convergence.

\vgap

We now present the convergence results for the stochastic APD method
applied to stochastic saddle-point problems with possibly unbounded feasible set $Z$.
It appears that the solution methods of these types of problems have not been well-studied in the literature.

\begin{theorem}
	\label{thmConvergenceUBS}
Let $\{z\ag\} = \{(x\ag, y\ag)\}$ be the iterates generated by Algorithm \ref{algACPD} 
with $V_X(x, x_t) = \|x - x_t\|^2/2$ and $V_Y(y, y_t) = \|y - y_t\|^2/2$. 
Assume that the parameters $\beta_t, \theta_t, \eta_t$ and $\tau_t$ 
in Algorithm \ref{algACPDS} satisfy \eqref{eqnIneqBetaGamma}, \eqref{eqnFlatStep} and
\eqnok{eqnIneqLGLKS} 
for all $t\ge 1$ and some $p,q\in(0,1)$, then there exists a perturbation vector $v_{t+1}$ such that 
				\begin{equation}
					\label{eqnErrUBDS}
					\E [\tilde g(z_{\tn}^{ag},v_{t+1})]\leq 
						\frac{1}{\beta_t\eta_t}\left(\frac{6-4p}{1-p}D^2 + \frac{5 - 3p}{2-2p}C^2 \right) =:\varepsilon\tn
				\end{equation}
				for any $t\geq 1$. Moreover, we have
				\begin{align}
					\label{eqnvboundS}
					& \E[\|v\tn\|]\leq  \frac {2\|\hat x - x_1\|}{\beta_t\eta_t} + \frac {2\|\hat y - y_1\|}{\beta_t\tau_t} + \sqrt{2D^2 + C^2} \left[\frac{2}{\beta_t\eta_t} + \frac{1}{\beta_t\tau_t}\sqrt{\frac{\tau_1}{\eta_1}}\left(\sqrt{\frac{1}{1-p}}+1\right) + \frac{2L_K}{\beta_t} \right],
				\end{align}
				where $(\hat x, \hat y)$ is a pair of solutions for problem \eqref{eqnDSPP}, $D$ is defined 
				in \eqnok{eqnD} and
				\begin{align}
					\label{eqnC}
					C := \sqrt{\sum_{i=1}^{t}\frac{\eta_i^2\sigma_x^2}{1-q} + \sum_{i=1}^{t}\frac{\eta_i\tau_i\sigma_y^2}{1-p}}.
				\end{align}
	
\end{theorem}

Below we specialize the results in Theorem~\ref{thmConvergenceUBS}
by choosing a set of parameters satisfying \eqref{eqnIneqBetaGamma}, \eqref{eqnFlatStep} and
\eqnok{eqnIneqLGLKS}.
\begin{corollary}
	\label{corStepsizeUBS}
	In Algorithm \ref{algACPDS}, if $N$ is given and the parameters are set to
	\begin{align}
		\label{eqnStepsizeUBS}
		\beta_t  =\frac{t+1}{2},\
		\theta_t  = \frac{t-1}{t},\
		\eta_t  = \frac{3t}{4\eta}, \text{ and }
		\tau_t = \frac{t}{\eta},
	\end{align}	
	where
	\begin{equation}
		\label{eqneta}
		\eta = 2L_G + 2L_K(N-1) + N\sqrt{N-1}\sigma / \tilde D \text{ for some }\tilde D>0,\ \sigma = \sqrt{\frac{9}{4}\sigma_x^2 + \sigma_y^2},
	\end{equation}
	then there exists $v_N$ that satisfies \eqref{eqnErrUBDS} with
		\begin{align}
			\varepsilon_N
			\leq&\  \frac{36L_GD^2}{N(N-1)} + \frac{36L_KD^2}{N} + \frac{\sigma D\left(18 D/\tilde D + 3\tilde D/D\right)}{\sqrt {N-1}} \label{sub1},
			\\
			\E[\|v_N\|] \leq&\ \frac{50L_GD}{N(N-1)} + \frac{L_KD(55 + 3\tilde D/D)}{N} + \frac{\sigma(6 + 25D/\tilde D)}{\sqrt{N-1}},\label{sub2}
		\end{align}
		where $D$ is defined in \eqref{eqnD}.
	\begin{proof}
		For the parameters in \eqref{eqnStepsizeUBS}, it is clear that \eqref{eqnIneqBetaGamma} and \eqref{eqnFlatStep} hold. Furthermore, let $p = 1/4$, $q= 3/4$, then for all $t=1,\ldots, N-1$, we have
		\begin{align*}
			\frac{q}{\eta_t}-\frac{L_G}{\beta_t}-\frac{L_K^2\tau_{t}}{p} = \frac{\eta}{t} - \frac{2L_G}{t+1} - \frac{4L_K^2t}{\eta} 
			\geq &\ \frac{2L_G + 2L_K(N-1)}{t} - \frac{2L_G}{t} - \frac{2L_K^2t}{L_K(N-1)} \geq 0,
		\end{align*}
		thus \eqref{eqnIneqLGLKS} holds. By Theorem \ref{thmConvergenceUBS}, we get \eqref{eqnErrUBDS} and \eqref{eqnvboundS}. Note that $\eta_t/\tau_t = 3/4$, and
		\begin{align*}
			 \frac{1}{\beta_{N-1}\eta_{N-1}}\|\hat x - x_1\| \leq \frac{1}{\beta_{N-1}\eta_{N-1}}D, \ \frac{1}{\beta_{N-1}\tau_{N-1}}\|\hat y - y_1\|\leq \frac{1}{\beta_{N-1}\eta_{N-1}}\cdot\frac{\eta_{N-1}}{\tau_{N-1}}\cdot\sqrt{\frac{4}{3}}D = \frac{\sqrt{3/4}D}{\beta_{N-1}\eta_{N-1}},\\
		\end{align*}
		so in \eqref{eqnErrUBDS} and \eqref{eqnvboundS} we have
		\begin{align}
			\label{tmpeps}
			\varepsilon\tn & \leq \frac{1}{\beta_{N-1}\eta_{N-1}}(\frac{20}{3}D^2 + \frac{17}{6}C^2),
			\\
			\label{tmpv}
			\ \E[\|v\tn\|] & \leq \frac{(2+\sqrt{3})D}{\beta_{N-1}\eta_{N-1}} + \frac{\sqrt{2D^2 + C^2}\left(3 + \sqrt{3/4} \right)}{\beta_{N-1}\eta_{N-1}} + \frac{2L_K\sqrt{2D^2 + C^2}}{\beta_{N-1}}.
		\end{align}
		By \eqref{eqnC} and the fact that $\sum_{i=1}^{N-1}i^2\leq N^2(N-1)/3$, we have
		\begin{align*}
			C = \sqrt{\sum_{i=1}^{N-1}\frac{9\sigma_x^2i^2}{4\eta^2} +  \sum_{i=1}^{N-1}\frac{\sigma_x^2i^2}{\eta^2}}\leq \sqrt{\frac{1}{3\eta^2}N^2(N-1)\left(\frac{9\sigma_x^2}{4} + {\sigma_x^2} \right)} = \frac{\sigma N\sqrt{N-1}}{\sqrt 3\eta}
		\end{align*}
		Applying the above bound to \eqref{tmpeps} and \eqref{tmpv}, and using the fact that $\sqrt{2D^2 + C^2}\leq \sqrt{2}D + C$, we obtain
		\begin{align*}
			&\eta = 2L_G + 2L_K(N-1) + N\sqrt{N-1}\sigma / \tilde D, \\
			\varepsilon_N \leq&\  \frac{8\eta}{3N(N-1)}\left(\frac{20}{3}D^2 + \frac{17\sigma^2N^2(N-1)}{18\eta^2}\right) = \frac{8}{3N(N-1)}\left(\frac{20}{3}\eta D^2 + \frac{17\sigma^2N^2(N-1)}{18\eta}\right)
			\\
			\leq&\  \frac{320L_GD^2}{9N(N-1)} + \frac{320L_K(N-1)D^2}{9N(N-1)} + \frac{160N\sqrt{N-1}\sigma D^2/\tilde D}{9N(N-1)} + \frac{68\sigma^2N^2(N-1)}{27N^2(N-1)^{3/2}\sigma/\tilde D}
			\\
			\leq&\  \frac{36L_GD^2}{N(N-1)} + \frac{36L_KD^2}{N} + \frac{\sigma D\left(18 D/\tilde D + 3\tilde D/D\right)}{\sqrt {N-1}},
			\\
			\E[\|v_N\|] \leq&\  \frac{1}{\beta_{N-1}\eta_{N-1}}\left(2D + \sqrt{3}D + 3\sqrt{2}D + \sqrt{6}D/2 + 3C + \sqrt{3}C/2\right) + \frac{2\sqrt{2}L_KD}{\beta_{N-1}} + \frac{2L_KC}{\beta_{N-1}}
			\\
			\leq &\ \frac{16L_G + 16L_K(N-1) + 8N\sqrt{N-1}\sigma/\tilde D}{3N(N-1)}\left(2 + \sqrt{3} + 3\sqrt{2} + \sqrt{6}/2\right)D
			\\
			&\ + \frac{8\sigma}{3\sqrt{N-1}}\left(\sqrt{3} + 1/2 \right) + \frac{4\sqrt{2}L_KD}{N} + \frac{4L_K\sigma N\sqrt{N-1}}{N\sqrt 3N\sqrt{N-1}\sigma/\tilde D}
			\\
			\leq&\ \frac{50L_GD}{N(N-1)} + \frac{L_KD(55 + 3\tilde D/D)}{N} + \frac{\sigma(6 + 25D/\tilde D)}{\sqrt{N-1}}.
		\end{align*}
		
	\end{proof}
\end{corollary}

Observe that the parameter settings in \eqnok{eqnStepsizeUBS}-\eqnok{eqneta}
are more complicated than the ones in \eqnok{eqnStepsizeUB}
for the deterministic unbounded case. In particular, for the stochastic unbounded case,
we need to choose a parameter $\tilde D$ which is not required for the deterministic case. 
Clearly, the optimal selection for $\tilde D$ minimizing the RHS of \eqnok{sub1} is given by $\sqrt{6} D$. 
Note however, that the value of $D$ will be very difficult to estimate for the unbounded case
and hence one often has to resort to a suboptimal selection for $\tilde D$.
For example, if $\tilde D = 1$, then the RHS of \eqnok{sub1} and \eqnok{sub2}
will become ${\cal O}(L_GD^2/N^2 + L_KD^2/N + \sigma D^2 / \sqrt{N})$
and ${\cal O}(L_G D/N^2 + L_KD/N + \sigma D / \sqrt{N})$, respectively.

\section{Convergence analysis} Our goal in this section is to prove the main results presented in 
Section 2 and 3, namely, Theorems~\ref{thmConvergenceBD}, \ref{thmConvergenceUB}, \ref{thmConvergenceBDS} and \ref{thmConvergenceUBS}.
\subsection{Convergence analysis for the deterministic APD algorithm}
\label{secProofDSPP}
In this section, we prove Theorems \ref{thmConvergenceBD} and \ref{thmConvergenceUB}
which, respectively, describe the convergence properties for the deterministic 
APD algorithm for the bounded and unbounded SPPs.

Before proving Theorem~\ref{thmConvergenceBD}, we first prove two technical 
results: Proposition~\ref{thmSimplifiedQ} shows some important
properties for the function $Q(\cdot, \cdot)$ in \eqnok{eqnQ}
and Lemma~\ref{lemQBoundGeneral} establishes a bound on
$Q(x^{ag}_{t}, z)$. 

\begin{proposition}
	\label{thmSimplifiedQ}
	Assume that $\beta_t\geq 1$ for all $t$. If $z\ag[1] = ( x\ag[1], y\ag[1])$ is generated by Algorithm \ref{algACPD}, then for all $z=(x, y)\in Z$,
	\begin{equation}
		\label{eqnSimplifiedQ}
		\begin{aligned}
				&\ \beta_tQ(z\ag[1], z) - (\beta_t-1)Q(z\ag, z)\\
				\leq &\ \langle\nabla G(x\md), x_{t+1}- x\rangle + \frac{L_G}{2\beta_t}\|x_{t+1}-x_t\|^2+\left[J(y_{t+1})-J( y)\right]  + \langle Kx_{t+1},  y\rangle - \langle K x, y_{t+1}\rangle.
		\end{aligned}
	\end{equation}	
	\begin{proof}
		By equations \eqref{eqnACPDxmd} and \eqref{eqnACPDxag}, $x\ag[1] - x\md = \beta_t^{-1}(x_{t+1}-x_t)$. Using this
		observation and the convexity of $G(\cdot)$, we have
		\begin{align}
			\label{eqnDerivationGxag}
			\begin{aligned}
				&\ \beta_tG(x\ag[1]) 
				\leq \  \beta_tG(x\md)+\beta_t\langle \nabla G(x\md), x\ag[1]-x\md\rangle + \frac{\beta_tL_G}{2}\|x\ag[1]-x\md\|^2
				\\
				\leq&\ \beta_tG(x\md) + \beta_t\langle\nabla G(x\md), x\ag[1]-x\md\rangle + \frac{L_G}{2\beta_t}\|x_{t+1}-x_t\|^2
				\\
				= &\ \beta_tG(x\md) + (\beta_t-1)\langle\nabla G(x\md), x\ag - x\md\rangle + \langle\nabla G(x\md), x_{t+1} - x\md\rangle + \frac{L_G}{2\beta_t}\|x_{t+1}-x_t\|^2
				\\
				= &\ (\beta_t-1)\left[G(x\md)+ \langle\nabla G(x\md), x\ag - x\md\rangle\right]  + \left[G(x\md)+\langle\nabla G(x\md), x_{t+1} - x\md\rangle\right] + \frac{L_G}{2\beta_t}\|x_{t+1}-x_t\|^2
				\\
				= &\ (\beta_t-1)\left[G(x\md)+ \langle\nabla G(x\md), x\ag - x\md\rangle\right]  +  \left[G(x\md)+\langle\nabla G(x\md),  x - x\md\rangle\right] + \langle\nabla G(x\md), x_{t+1} -  x\rangle \\
				& + \frac{L_G}{2\beta_t}\|x_{t+1}-x_t\|^2
				\\
				\leq&\ (\beta_t-1)G(x\ag) + G( x) +  \langle\nabla G(x\md), x_{t+1} -  x\rangle + \frac{L_G}{2\beta_t}\|x_{t+1}-x_t\|^2.
			\end{aligned}
		\end{align}	
Moreover, by \eqref{eqnACPDyag} and the convexity of $J(\cdot)$, we have
		\begin{align*}
			\begin{aligned}
			\beta_tJ(y\ag[1]) - \beta_tJ( y) &\ \leq (\beta_t-1)J(y\ag) + J(y_{t+1}) - \beta_tJ( y)\\
			&\ = \ (\beta_t-1)\left[J(y\ag)-J( y)\right] + J(y_{t+1})-J( y).
			\end{aligned}
		\end{align*}
By \eqref{eqnQ}, \eqref{eqnACPDxag}, \eqref{eqnACPDyag} and the above two inequalities above, we obtain
		\begin{align*}
			&\ \beta_tQ(z\ag[1], z) - (\beta_t-1)Q(z\ag, z)\\
			= &\ \beta_t\left\{\left[G(x\ag[1])+\langle Kx\ag[1],  y\rangle-J( y)\right] - \left[G( x)+\langle K x, y\ag[1]\rangle-J(y\ag[1])\right]\right\}\\
			& - (\beta_t-1)\left\{\left[G(x\ag) + \langle Kx\ag,   y\rangle-J( y)\right]  -\left[G( x)+\langle K x, y\ag\rangle-J(y\ag)\right]\right\}\\
			= &\ \beta_tG(x\ag[1])- (\beta_t-1)G(x\ag)-G( x) + \beta_t\left[J(y\ag[1])-J( y)\right]
			\\
			& - (\beta_t-1)\left[J(y\ag)-J( y)\right] + \langle K(\beta_tx\ag[1]-(\beta_t-1)x\ag),  y\rangle - \langle K x, \beta_ty\ag[1]-(\beta_t-1)y\ag\rangle\\
			\leq &\ \langle\nabla G(x\md), x_{t+1}- x\rangle + \frac{L_G}{2\beta_t}\|x_{t+1}-x_t\|^2+J(y_{t+1})-J( y) + \langle Kx_{t+1},  y\rangle - \langle K x, y_{t+1}\rangle.
		\end{align*}
	\end{proof}
\end{proposition}

Lemma~\ref{lemQBoundGeneral} establishes a bound for $Q(z\ag[1], z)$ for all $z \in Z$,
which will be used in the proof of both Theorems~\ref{thmConvergenceBD} and \ref{thmConvergenceUB}. 

\begin{lemma}
	\label{lemQBoundGeneral}
	Let $z\ag[1] = (x\ag[1], y\ag[1])$ be the iterates generated by Algorithm \ref{algACPD}. 
	Assume that the parameters $\beta_t, \theta_t, \eta_t$, and $\tau_t$ satisfy 
	\eqref{eqnIneqBetaGamma}, \eqref{eqnIncStep} and \eqref{eqnIneqLGLK}. 
	Then, for any $z\in Z$, we have
	\beq \label{eqnQBoundGeneral}
	\beta_t\gamma_t Q(z^{ag}_{t+1}, z) \le {\cal B}_t(z,z_{[t]}) +\gamma_{t}\langle K(x_{t+1}-x_{t}),  y-y_{t+1}\rangle
	- \gamma_{t}\left(\frac{\alpha_X}{2\eta_{t}}-\frac{L_G}{2\beta_{t}} \right)\|x_{t+1}-x_{t}\|^2,
	\eeq
	where $\gamma_t$ is defined in \eqref{eqnGamma}, $z_{[t]} := \{(x_i, y_i)\}_{i=1}^{t+1}$ and
		\begin{align}
		\label{def_CB}
		\begin{aligned}
			{\cal B}_t(z, z_{[t]}) &:= \sum_{i=1}^t
			\left\{\frac{\gamma_i}{\eta_i}\left[V_X(x, x_i) - V_X(x, x_{i+1}) \right]
+ \frac{\gamma_i}{\tau_i}\left[V_Y(y, y_i) - V_Y(y, y_{i+1}) \right] \right\}.
		\end{aligned}
		\end{align}	

	\begin{proof}
	First of all, we explore the optimality conditions in iterations \eqref{eqnACPDyt1} and \eqref{eqnACPDxt1}. Apply Lemma 2 in \cite{GhaLan10-1b} to \eqref{eqnACPDyt1}, we have
	\begin{align}
		\label{eqnACPDyt1Condition}
		\begin{aligned}
			&\ \langle-K\bar x_t, y_{t+1}- y\rangle + J(y_{t+1})-J( y) \leq \frac 1{\tau_t}V_Y( y,y_t) - \frac 1{\tau_t}V_Y(y_{t+1},y_t) - \frac {1}{\tau_t}V_Y( y, y_{t+1})\\
			\leq &\ \frac 1{\tau_t}V_Y( y,y_t) - \frac {\alpha_Y}{2\tau_t}\|y_{t+1}-y_t\|^2 - \frac {1}{\tau_t}V_Y( y, y_{t+1}),
		\end{aligned}
	\end{align}
	where the last inequality follows from the fact that, by the strong convexity of $d_Y(\cdot)$ and \eqnok{eqnVy},
	\begin{align}
		\label{eqnVNorm}
		V_Y(y_1, y_2)\geq \frac{\alpha_Y}2\|y_1-y_2\|^2, \text{ for all }y_1, y_2\in \cal Y.
	\end{align}
	Similarly, from \eqref{eqnACPDxt1} we can derive that
	\begin{align}
		\label{eqnACPDxt1Condition}
		\begin{aligned}
			& \langle \nabla G(x\md), x_{t+1}- x\rangle + \langle x_{t+1}- x, K^Ty\tn\rangle \leq \frac 1{\eta_t}V_X( x,x_t) - \frac {\alpha_X}{2\eta_t}\|x_{t+1} - x_t\|^2 - \frac {1}{\eta_t}V_X( x,x_{t+1}).\\
		\end{aligned}
	\end{align}
Our next step is to establish a crucial recursion of Algorithm \ref{algACPD}. 
	It follows from	\eqref{eqnSimplifiedQ}, \eqref{eqnACPDyt1Condition} and \eqref{eqnACPDxt1Condition} that
	\begin{align}
	\label{temp_label}
		\begin{aligned}
				&\ \beta_t Q(z\ag[1], z) - (\beta_t-1) Q(z\ag, z)
				\\
				\leq &\ \langle\nabla G(x\md), x_{t+1}- x\rangle + \frac{L_G}{2\beta_t}\|x_{t+1}-x_t\|^2+\left[J(y_{t+1})-J( y)\right]  + \langle Kx_{t+1},  y\rangle - \langle K x, y_{t+1}\rangle
				\\
				\leq &\ \frac {1}{\eta_t}V_X( x,x_t) - \frac {1}{\eta_t}V( x,x_{t+1})-\left( \frac {\alpha_X}{2\eta_t}- \frac{L_G}{2\beta_t}\right)\|x_{t+1}-x_t\|^2\\
				&\ + \frac {1}{\tau_t}V_Y( y,y_t) - \frac {1}{\tau_t}V( y,y_{t+1})- \frac {\alpha_Y}{2\tau_t}\|y_{t+1}-y_t\|^2  \\
				&\ -\langle x_{t+1}- x, K^Ty_{t+1}\rangle
				+\langle K\bar x_t,y_{t+1}- y\rangle + \langle Kx_{t+1},  y\rangle - \langle K x, y_{t+1}\rangle.
			\end{aligned}
		\end{align}
Also observe that by \eqref{eqnACPDxbar}, we have
		\begin{align*}
				&\ -\langle x_{t+1}- x, K^Ty_{t+1}\rangle
				 + \langle K\bar x_t,y_{t+1}- y\rangle + \langle Kx_{t+1},  y\rangle - \langle K x, y_{t+1}\rangle\\
 				= &\ \langle K(x_{t+1}-x_t),  y-y_{t+1}\rangle-\theta_{t}\langle K(x_t-x_{t-1}), y-y_{t+1}\rangle\\
 				= &\ \langle K(x_{t+1}-x_t),  y-y_{t+1}\rangle-\theta_{t}\langle K(x_t-x_{t-1}), y-y_{t}\rangle -\theta_{t}\langle K(x_t-x_{t-1}),y_{t}-y_{t+1}\rangle.
		\end{align*}
Multiplying both sides of \eqnok{temp_label} by $\gamma_t$,
using the above identity and the fact that $\gamma_t \theta_t = \gamma_{t-1}$ due to \eqnok{eqnGamma},
we obtain
		\begin{equation}
			\label{eqnACPDRecursion}
			\begin{aligned}
				&\ \beta_t\gamma_tQ(z\ag[1], z) - (\beta_t-1)\gamma_tQ(z\ag, z)\\
				\leq &\ \frac {\gamma_t}{\eta_t}V_X( x, x_t) - \frac {\gamma_t}{\eta_t}V_X( x,x_{t+1}) + \frac {\gamma_t}{\tau_t}V_Y( y, y_t) - \frac {\gamma_t}{\tau_t}V_Y( y, y_{t+1}) \\
				& +\gamma_t\langle K(x_{t+1}-x_{t}),  y-y_{t+1}\rangle - \gamma_{t-1}\langle K(x_t-x_{t-1}),  y-y_t\rangle\\
				& -\gamma_t\left(\frac{\alpha_X}{2\eta_t} - \frac{L_G}{2\beta_t}\right)\|x_{t+1}-x_t\|^2 - \frac {\alpha_Y\gamma_t}{2\tau_t}\|y_{t+1}-y_t\|^2  -\gamma_{t-1}\langle K(x_t-x_{t-1}),y_{t}-y_{t+1}\rangle.
			\end{aligned}
		\end{equation}		
		Now, applying Cauchy-Schwartz inequality to the last term in \eqref{eqnACPDRecursion}, using the notation $L_K=\|K\|$ and noticing that $\gamma_{t-1}/\gamma_{t} = \theta_t \leq \min\{\eta_{t-1}/\eta_{t},  \tau_{t-1}/\tau_{t} \}$ from \eqref{eqnIncStep}, we have
			\begin{align}
				\label{eqnCS}
				\begin{aligned}
					&\ -\gamma_{t-1}\langle K(x_t-x_{t-1}),y_{t}-y_{t+1}\rangle
					\leq \gamma_{t-1}\|K(x_t-x_{t-1})\|_*\|y_t-y_{t+1}\|\\
					\leq &\ L_K\gamma_{t-1}\|x_t-x_{t-1}\| \, \|y_t-y_{t+1}\|
					\leq  \frac{L_K^2\gamma_{t-1}^2\tau_{t}}{2\alpha_Y\gamma_t}\|x_t-x_{t-1}\|^2+\frac{\alpha_Y\gamma_t}{2\tau_{t}}\|y_{t}-y_{t+1}\|^2
					\\
					\leq &\ 
					\frac{L_K^2\gamma_{t-1}\tau_{t-1}}{2\alpha_Y}\|x_t-x_{t-1}\|^2+\frac{\alpha_Y\gamma_t}{2\tau_{t}}\|y_{t}-y_{t+1}\|^2.
				\end{aligned}	
			\end{align}
Noting that $\theta_{t+1} = \gamma_t/\gamma\tn$, so by \eqref{eqnIneqBetaGamma} we have $(\beta\tn -1)\gamma\tn = \beta_t\gamma_t$.
Combining the above two relations with inequality \eqref{eqnACPDRecursion}, we get
		the following recursion for Algorithm~\ref{algACPD}.
	\begin{equation*}
		\begin{aligned}
			&\ (\beta_{t+1}-1)\gamma_{t+1}Q(z\ag[1], z) - (\beta_t-1)\gamma_tQ(z\ag, z)	
			= \ \beta_t\gamma_tQ(z\ag[1], z) - (\beta_t-1)\gamma_tQ(z\ag, z)\\
			\leq &\ \frac {\gamma_t}{\eta_t}V_X( x, x_t) - \frac {\gamma_t}{\eta_t}V_X( x,x_{t+1}) + \frac {\gamma_t}{\tau_t}V_Y( y, y_t) - \frac {\gamma_t}{\tau_t}V_Y( y, y_{t+1}) \\
			& +\gamma_t\langle K(x_{t+1}-x_{t}),  y-y_{t+1}\rangle - \gamma_{t-1}\langle K(x_t-x_{t-1}),  y-y_t\rangle\\
			& -\gamma_t\left(\frac{\alpha_X}{2\eta_t} -  \frac{L_G}{2\beta_t}\right)\|x_{t+1}-x_t\|^2 + \frac{L_K^2\gamma_{t-1}\tau_{t-1}}{2\alpha_Y}\|x_t-x_{t-1}\|^2, \forall t \ge 1.
		\end{aligned}
	\end{equation*}	
Applying the above inequality inductively and assuming that $x_0 = x_1$, we conclude that
	\begin{align*}
			\begin{aligned}
			&	(\beta_{t+1}-1)\gamma_{t+1} Q(z^{ag}_{t+1}, z) - (\beta_1 - 1)\gamma_1Q(z^{ag}_1, z) \le  {\cal B}_t(z, z_{[t]}) +\gamma_{t}\langle K(x_{t+1}-x_{t}),  y-y_{t+1}\rangle				  
				\\		
				& -\gamma_{t}\left(\frac{\alpha_X}{2\eta_{t}}-\frac{L_G}{2\beta_{t}} \right)\|x_{t+1}-x_{t}\|^2
				- \sum_{i=1}^{t-1}\gamma_i\left(\frac{\alpha_X}{2\eta_i} - \frac{L_G}{2\beta_i} - \frac{L_K^2\tau_i}{2\alpha_Y} \right)\|x_{i+1} - x_i\|^2,
			\end{aligned}
			\end{align*}
			which, in view of \eqref{eqnIneqLGLK}
			and the facts that $\beta_1 = 1$ and $(\beta\tn -1)\gamma\tn = \beta_t\gamma_t$ by \eqref{eqnIneqBetaGamma},
			implies \eqref{eqnQBoundGeneral}.
		\end{proof}
	\end{lemma}

\vgap

	We are now ready to prove Theorem \ref{thmConvergenceBD}, which follows as an immediate consequence of Lemma \ref{lemQBoundGeneral}.

\noindent	{\bf Proof of Theorem~\ref{thmConvergenceBD}}.
Let ${\cal B}_t(z)$ be defined in \eqnok{def_CB}. 
First note that by the definition of $\gamma_t$ in \eqnok{eqnGamma} and relation \eqnok{eqnIncStep},
we have $\theta_t = \gamma_{t-1}/\gamma_t \le \eta_{t-1}/\eta_t$ and hence $\gamma_{t-1}/\eta_{t-1}
\le \gamma_t / \eta_t$. Using this observation and \eqnok{boundness}, we conclude that
\begin{align}
\label{bound_CB}
			\begin{aligned}
{\cal B}_t(z, z_{[t]}) &=
\frac {\gamma_1}{\eta_1}V_X(x, x_1) - \sum_{i=1}^{t-1}\left(\frac {\gamma_i}{\eta_i} - \frac {\gamma_{i+1}}{\eta_{i+1}}\right)V_X(x,x_{i+1}) - 
\frac{\gamma_{t}}{\eta_{t}}V_X( x, x_{t+1})
			\\
			& \ \ +\frac{\gamma_1}{\tau_1}V_Y(y, y_1) - \sum_{i=1}^{t-1}\left(\frac {\gamma_i}{\tau_i} - \frac {\gamma_{i+1}}{\tau_{i+1}}\right)V_Y( y,y_{i+1}) - \frac {\gamma_{t}}{\tau_{t}}V_Y( y, y_{t+1})\\
& \le \frac {\gamma_1}{\eta_1}\Omega_X^2 - \sum_{i=1}^{t-1}\left(\frac {\gamma_i}{\eta_i} - \frac {\gamma_{i+1}}{\eta_{i+1}}\right)\Omega_X^2 - \frac{\gamma_{t}}{\eta_{t}}V_X( x, x_{t+1})\\
		&	\ \ +
			\frac{\gamma_1}{\tau_1}\Omega_Y^2 - \sum_{i=1}^{t-1}\left(\frac {\gamma_i}{\tau_i} - \frac {\gamma_{i+1}}{\tau_{i+1}}\right)\Omega_Y^2
			- \frac {\gamma_{t}}{\tau_{t}}V_Y( y, y_{t+1})\\
			& = \frac {\gamma_t}{\eta_t}\Omega_X^2 - \frac{\gamma_{t}}{\eta_{t}}V_X( x, x_{t+1})
			+\frac{\gamma_t}{\tau_t}\Omega_Y^2 - \frac {\gamma_{t}}{\tau_{t}} V_Y( y, y_{t+1}).
			\end{aligned}
\end{align}
Now applying Cauchy-Schwartz inequality
	to the inner product term in \eqref{eqnQBoundGeneral}, we get
		\begin{align}
		\label{bound_inner}
			\begin{aligned}
			&\ \gamma_{t}\langle K(x_{t+1}-x_{t}),  y-y_{t+1}\rangle
			\leq L_K\gamma_{t}\|x_{t+1}-x_{t}\|\| y-y_{t+1}\|
			\leq  \frac{L_K^2\gamma_{t}\tau_{t}}{2\alpha_Y}\|x_{t+1}-x_{t}\|^2 + \frac{\alpha_Y\gamma_{t}}{2\tau_{t}}\| y-y_{t+1}\|^2.
			\end{aligned}
		\end{align}	
Using the above two relations, \eqref{eqnIneqLGLK}, \eqref{eqnQBoundGeneral}
and \eqnok{eqnVNorm}, we have
		\begin{align}
		\begin{aligned}
			&\ \beta_t\gamma_t Q(z^{ag}_{t+1}, z) \leq \ \frac {\gamma_t}{\eta_t}\Omega_X^2 - \frac{\gamma_{t}}{\eta_{t}}V_X( x, x_{t+1})
			+\frac{\gamma_t}{\tau_t}\Omega_Y^2 - \frac {\gamma_{t}}{\tau_{t}}\left(V_Y( y, y_{t+1}) - \frac{\alpha_Y}{2}\| y-y_{t+1}\|^2\right) 
			\\
			&  -\gamma_{t}\left(\frac{\alpha_X}{2\eta_{t}}-\frac{L_G}{2\beta_{t}}-\frac{L_K^2\tau_{t}}{2\alpha_Y} \right)\|x_{t+1}-x_{t}\|^2 
			\le   \frac {\gamma_t}{\eta_t}\Omega_X^2 + \frac{\gamma_t}{\tau_t}\Omega_Y^2, \ \ \forall z \in Z,
		\end{aligned}
		\end{align}		
which together with \eqnok{def_gap}, then clearly imply \eqref{eqnErrBDD}.
	\endproof

\vgap

Our goal in the remaining part of this subsection is to prove Theorem~\ref{thmConvergenceUB},
which summarizes the convergence properties of Algorithm \ref{algACPD} when $X$ or $Y$ is unbounded.
We will first prove a technical result which specializes the results in Lemma~\ref{lemQBoundGeneral}
for the case when \eqref{eqnIneqBetaGamma}, \eqref{eqnFlatStep} and \eqref{eqnIneqLGLKp} hold.
	
	\begin{lemma}
		\label{lemEucl}
		Let $\hat z = (\hat x, \hat y)\in Z$ be a saddle point of \eqref{eqnDSPP}. If $V_X(x, x_t)=\|x - x_t\|^2/2$ 
		and $V_Y(y, y_t)=\|y - y_t\|^2/2$ in Algorithm \ref{algACPD}, and the parameters $\beta_t, \theta_t, \eta_t$
		and $\tau_t$ satisfy \eqref{eqnIneqBetaGamma}, \eqref{eqnFlatStep} and \eqref{eqnIneqLGLKp}, then
		\begin{enumerate}
			\renewcommand{\theenumi}{(\alph{enumi})}
			\itemeqn 
			\begin{equation}
				\label{eqnIneqSadPt}
				\|\hat x - x\tn\|^2 + \frac{\eta_t(1-p)}{\tau_t}\|\hat y - y\tn\|^2 \leq \|\hat x - x_1\|^2 + \frac{\eta_t}{\tau_t}\|\hat y - y_1\|^2, \text{ for all }t\geq 1.
			\end{equation}
			
			\itemeqn
			\begin{equation}
				\label{eqnModifiedQ}
				\tilde g(z_{t+1}^{ag}, v_{t+1}) \leq \frac 1{2\beta_t\eta_t}\|x\ag[1] - x_1\|^2 + \frac 1{2\beta_t\tau_t}\|y\ag[1] - y_1\|^2 =:\delta_{t+1}, \text{ for all }t\geq 1,
			\end{equation}
			where $\tilde g(\cdot, \cdot)$ is defined in \eqnok{def_gapt} and
			\begin{align}
				\label{eqnv}
				v_{t+1} &=\left(\frac 1{\beta_t\eta_t}(x_1 - x\tn) , \frac 1{\beta_t\tau_t}(y_1 - y\tn) + \frac 1{\beta_t} K(x\tn - x_t) \right).
			\end{align}
		\end{enumerate}
	\begin{proof}
		It is easy to check that the conditions in Lemma \ref{lemQBoundGeneral} are satisfied. By \eqref{eqnFlatStep},  \eqref{eqnQBoundGeneral} in Lemma \ref{lemQBoundGeneral} becomes
		\begin{align}
			\label{eqnQBoundEucl}
		\begin{aligned}
			\beta_t Q(z^{ag}_{t+1}, z) 		
			\leq &\ \frac {1}{2\eta_t}\|x - x_1\|^2  - \frac{1}{2\eta_{t}}\| x- x_{t+1}\|^2 +\frac{1}{2\tau_t}\|y- y_1\|^2  - \frac {1}{2\tau_{t}}\|y- y_{t+1}\|^2 
			\\
			& +\langle K(x_{t+1}-x_{t}),  y-y_{t+1}\rangle -\left(\frac{1}{2\eta_{t}}-\frac{L_G}{2\beta_{t}} \right)\|x_{t+1}-x_{t}\|^2.
		\end{aligned}
		\end{align}	
		
		To prove \eqref{eqnIneqSadPt}, observe that
		\begin{align}
			\label{eqnLastProdp}
			\begin{aligned}
			&\ \langle K(x_{t+1}-x_{t}),  y-y_{t+1}\rangle
			\leq \ \frac{L_K^2\tau_{t}}{2p}\|x_{t+1}-x_{t}\|^2 + \frac{p}{2\tau_{t}}\| y-y_{t+1}\|^2
			\end{aligned}
		\end{align}			
		where $p$ is the constant in \eqref{eqnIneqLGLKp}. 
		By \eqref{eqnIneqLGLKp} and the above two inequalities, we get
		\begin{align*}
		\begin{aligned}
			\beta_t Q(z^{ag}_{t+1}, z) 		
			\leq &\ \frac {1}{2\eta_t}\|x - x_1\|^2  - \frac{1}{2\eta_{t}}\| x- x_{t+1}\|^2 +\frac{1}{2\tau_t}\|y- y_1\|^2  - \frac {1-p}{2\tau_{t}}\|y- y_{t+1}\|^2.
		\end{aligned}
		\end{align*}		
		Letting $z=\hat z$ in the above, and using the fact that $Q(z\ag[1], \hat z)\geq 0$, we obtain \eqref{eqnIneqSadPt}.
		
		Now we prove \eqref{eqnModifiedQ}. Noting that
		\begin{align}
			\label{eqnx1toxag}
			\begin{aligned}
			& \|x - x_1\|^2 - \|x - x\tn\|^2  = 2\langle x\tn - x_1, x\rangle + \|x_1\|^2 - \|x\tn\|^2\\
			= & 2\langle x\tn - x_1, x - x\ag[1]\rangle + 2\langle x\tn - x_1, x\ag[1]\rangle + \|x_1\|^2 - \|x\tn\|^2\\
			= & 2\langle x\tn - x_1, x - x\ag[1]\rangle + \|x\ag[1] - x_1\|^2 - \|x\ag[1] - x\tn\|^2,
			\end{aligned}
		\end{align}		
		we conclude from \eqref{eqnIneqLGLKp} and \eqref{eqnQBoundEucl} that for any $z\in Z$,
		\begin{align*}
			&\ \beta_tQ(z\ag[1], z) - \langle K(x\tn - x_t), y\ag[1] - y\rangle - \frac 1{\eta_t}\langle x_1 - x\tn, x\ag[1] - x\rangle - \frac 1{\tau_t}\langle y_1 - y\tn, y\ag[1] - y\rangle 
			\\
			\leq & \ \frac 1{2\eta_t}\left(\|x\ag[1] - x_1\|^2 - \|x\ag[1] - x\tn\|^2 \right) + \frac 1{2\tau_t}\left(\|y\ag[1] - y_1\|^2 - \|y\ag[1] - y\tn\|^2 \right)
			\\
			&\ +\langle K(x\tn - x_t), y\ag[1] - y\tn\rangle - \left(\frac{1}{2\eta_{t}}-\frac{L_G}{2\beta_{t}} \right)\|x_{t+1}-x_{t}\|^2
			\\
			\leq & \ \frac 1{2\eta_t}\left(\|x\ag[1] - x_1\|^2 - \|x\ag[1] - x\tn\|^2 \right) + \frac 1{2\tau_t}\left(\|y\ag[1] - y_1\|^2 - \|y\ag[1] - y\tn\|^2 \right)
			\\
			&\ \frac{p}{2\tau_t}\|y\ag[1] - y\tn\|^2 -  \left(\frac{1}{2\eta_{t}}-\frac{L_G}{2\beta_{t}} - \frac{L_K^2\tau_t}{2p} \right)\|x_{t+1}-x_{t}\|^2
			\\
			\leq & \ \frac 1{2\eta_t}\|x\ag[1] - x_1\|^2 + \frac 1{2\tau_t}\|y\ag[1] - y_1\|^2.
		\end{align*}
The result in \eqref{eqnModifiedQ} and \eqref{eqnv} immediately follows from the above inequality
and \eqnok{def_gapt}.
		
	\end{proof}
	
	\end{lemma}
	
We are now ready to prove	Theorem \ref{thmConvergenceUB}.
	
\noindent{\bf Proof of Theorem \ref{thmConvergenceUB}}. 
		We have established the expression of $v\tn$ and $\delta\tn$ in Lemma \ref{lemEucl}. It
		suffices to estimate the bound on $\|v\tn\|$ and $\delta\tn$.
		It follows from the definition of $D$, \eqref{eqnFlatStep} and \eqref{eqnIneqSadPt} that for all $t\ge 1$,
	\[
			\|\hat x - x\tn\| \leq D \textrm{ and \ } \|\hat y - y\tn\| \leq D\sqrt{\frac{\tau_1}{\eta_1(1-p)}}.
		\]
		Now by \eqref{eqnv}, we have
		\begin{align*}
			\|v_{t+1}\| & \leq \frac 1{\beta_t\eta_t}\|x_1 - x\tn\| + \frac 1{\beta_t\tau_t}\|y_1 - y\tn\| + \frac {L_K}{\beta_t} \|x\tn - x_t\|
			\\
			& \leq \frac 1{\beta_t\eta_t}\left(\|\hat x - x_1\| + \|\hat x - x\tn\|\right) + \frac 1{\beta_t\tau_t}\left(\|\hat y - y_1\| + \|\hat y  - y\tn\|\right) + \frac {L_K}{\beta_t} \left(\|\hat x - x\tn\| + \|\hat x - x_t\|\right)
			\\
			& \leq \frac 1{\beta_t\eta_t}\left(\|\hat x - x_1\| + D\right) + \frac 1{\beta_t\tau_t}\left(\|\hat y - y_1\| + D\sqrt{\frac{\tau_1}{\eta_1(1-p)}}\right) + \frac {2L_K}{\beta_t} D\\
			& = \frac 1{\beta_t\eta_t}\|\hat x - x_1\| + \frac 1{\beta_t\tau_t}\|\hat y - y_1\| + D \left[\frac{1}{\beta_t\eta_t}\left(1 + \sqrt{\frac{\eta_1}{\tau_1(1-p)}}\right) + \frac{2L_K}{\beta_t} \right].
		\end{align*}
To estimate the bound of $\delta\tn$, consider the sequence $\{\gamma_t \}$ defined in \eqnok{eqnGamma}. 
Using the fact that $(\beta\tn - 1)\gamma\tn = \beta_t\gamma_t$ due to \eqnok{eqnIneqBetaGamma} and
\eqnok{eqnGamma}, and applying \eqref{eqnACPDxag} and \eqref{eqnACPDyag} inductively, we have
		\beq \label{eqnConvex_comb}
		x\ag[1] = \frac{1}{\beta_t\gamma_t}\sum_{i = 1}^{t}\gamma_ix\tn,\ y\ag[1] = \frac{1}{\beta_t\gamma_t}\sum_{i = 1}^{t}\gamma_iy\tn \ \textrm{and }\ \frac{1}{\beta_t\gamma_t}\sum_{i=1}^{t}\gamma_i = 1. \eeq
Thus $x\ag[1]$ and $y\ag[1]$ are convex combinations of sequences $\{x_{i+1} \}_{i=1}^t$ and $\{y_{i+1} \}_{i=1}^t$ . 
Using these relations and \eqref{eqnIneqSadPt}, we have
		\begin{align*}
			\delta\tn & = \frac 1{2\beta_t\eta_t}\|x\ag[1] - x_1\|^2 + \frac 1{2\beta_t\tau_t}\|y\ag[1] - y_1\|^2
			\\ 
			& \leq \frac{1}{\beta_t\eta_t}\left(\|\hat x - x\ag[1]\|^2 + \|\hat x - x_1\|^2  \right) + \frac{1}{\beta_t\tau_t}\left(\|\hat y - y\ag[1]\|^2 + \|\hat y - y_1\|^2  \right)
			\\
			& = \frac{1}{\beta_t\eta_t}\left(D^2 + \|\hat x - x\ag[1]\|^2 + \frac{\eta_t(1-p)}{\tau_t}\|\hat y - y\ag[1]\|^2 + \frac{\eta_tp}{\tau_t}\|\hat y - y\ag[1]\|^2\right)
			\\
			& \leq \frac{1}{\beta_t\eta_t}\left[D^2 + \frac{1}{\beta_t\gamma_t}\sum_{i=1}^{t}\gamma_i\left(\|\hat x - x_{i+1}\|^2 + \frac{\eta_t(1-p)}{\tau_t}\|\hat y - y_{i+1}\|^2  + \frac{\eta_tp}{\tau_t}\|\hat y - y_{i+1}\|^2\right)\right]
			\\
			& \leq \frac{1}{\beta_t\eta_t}\left[D^2 + \frac{1}{\beta_t\gamma_t}\sum_{i=1}^{t}\gamma_i\left(D^2 + \frac{\eta_tp}{\tau_t}\cdot\frac{\tau_1}{\eta_1(1-p)}D^2 \right)\right] = \frac{(2-p)D^2}{\beta_t\eta_t(1-p)}.
		\end{align*}

	\endproof

\subsection{Convergence analysis for the stochastic APD algorithm}

In this subsection, we prove Theorems \ref{thmConvergenceBDS} and  \ref{thmConvergenceUBS}
which descirbe the convergence properties of the stochastic APD algorithm
presented in Section 3.

Let $\gi$, $\pi$ and $\ri$ be the output from the $\SO$ at the $t$-th iteration
of Algorithm~\ref{algACPDS}. Throughout this subsection, we denote 
\[
\Delta_{x,G}^t:=\gi - \nabla G(x_t^{md}), \ \Delta_{x,K}^t:=\ri - K^T y_{t+1}, \  \Delta_{y}^t:=- \pi + K\bar x_t, 
\]
\[
\Delta_x^t := \Delta_{x,G}^t+\Delta_{x,K}^t \ \ \mbox{and} \ \ \Delta^t := (\Delta_x^t, \Delta_y^t).
\]
Moreover,
for a given $z= (x, y) \in Z$, let us denote $\|z\|^2 = \|x\|^2 + \|y\|^2$
and its associate dual norm for $\Delta = (\Delta_x, \Delta_y)$
by $\|\Delta\|_*^2 = \|\Delta_x\|_*^2 + \|\Delta_y\|_*^2$.
We also define the Bregman divergence $V(z, \tilde z)
:= V_X(x, \tilde x) + V_Y(y, \tilde y)$ for $z = (x, y)$
and $\tilde z = (\tilde x, \tilde y)$.

\vgap

Before proving Theorem~\ref{thmConvergenceBDS}, we first estimate a bound on $Q(z\ag[1], z)$ for all $z \in Z$. 
This result is analogous to Lemma \ref{lemQBoundGeneral} for the deterministic APD method.

\begin{lemma}
	\label{lemQBoundGeneralS}
Let $z^{ag}_t = (x^{ag}_t, y^{ag}_t)$ be the iterates generated by Algorithm \ref{algACPDS}.
Assume that the parameters $\beta_t, \theta_t, \eta_t$ and $\tau_t$ satisfy 
\eqref{eqnIneqBetaGamma}, \eqref{eqnIncStep} and \eqref{eqnIneqLGLKS}. Then, for any $z\in Z$, we have
	\begin{align}
		\label{eqnQBoundGeneralS}
		\begin{aligned}
			&\ \beta_t\gamma_t Q(z^{ag}_{t+1}, z) \leq {\cal B}_t(z, z_{[t]}) +\gamma_{t}\langle K(x_{t+1}-x_{t}),  y-y_{t+1}\rangle -\gamma_{t}\left(\frac{q\alpha_X}{2\eta_{t}}-\frac{L_G}{2\beta_{t}} \right)\|x_{t+1}-x_{t}\|^2 + \sum_{i=1}^{t}\Lambda_i(z),
		\end{aligned}
	\end{align}
	where $\gamma_t$ and ${\cal B}_t(z, z_{[t]})$, respectively, are defined in \eqnok{eqnGamma} and \eqnok{def_CB}, 
	$z_{[t]} = \{(x_i,y_i)\}_{i=1}^{t+1}$ and
		\begin{align}
			\label{eqnLambda}
			\Lambda_i(z) & := -\frac{(1-q)\alpha_X\gamma_i}{2\eta_i}\|x_{i+1}-x_i\|^2 - \frac {(1-p)\alpha_Y\gamma_i}{2\tau_i}\|y_{i+1}-y_i\|^2
			-\gamma_i\langle\Delta^i, z_{i+1}-  z\rangle.
		\end{align}

\begin{proof}
	Similar to \eqref{eqnACPDyt1Condition} and \eqref{eqnACPDxt1Condition}, we conclude from the optimality conditions of  \eqref{eqnACPDSyt1} and \eqref{eqnACPDSxt1} that
			\begin{align*}
				\begin{aligned}
					 \langle-\pt, y_{t+1}-y\rangle + J(y_{t+1})-J(\hat y) 				\leq &\ \frac 1{\tau_t}V_Y(\hat y,y_t) - \frac {\alpha_Y}{2\tau_t}\|y_{t+1}-y_t\|^2 - \frac {1}{\tau_t}V_Y(\hat y, y_{t+1}),\\
					 \langle \gt, x_{t+1}-x\rangle + \langle x_{t+1}-x, \rt\rangle
					\leq &\ \frac 1{\eta_t}V_X(x,x_t) - \frac {\alpha_X}{2\eta_t}\|x_{t+1}-x_t\|^2 - \frac {1}{\eta_t}V_X(x,x_{t+1}).
				\end{aligned}
			\end{align*}

		Now we establish an important recursion for Algorithm \ref{algACPDS}. Observing that Proposition \ref{thmSimplifiedQ} also holds
		for Algorithm \ref{algACPDS}, and applying the above two inequalities to \eqref{eqnSimplifiedQ} in Proposition \ref{thmSimplifiedQ}, similar to  \eqref{eqnACPDRecursion}, we have
		\begin{equation}
			\label{eqnACPDSRecursion}
			\begin{aligned}
				&\ \beta_t\gamma_tQ(z\ag[1], z) - (\beta_t-1)\gamma_tQ(z\ag, z)\\
				\leq &\ \frac {\gamma_t}{\eta_t}V_X(x, x_t) - \frac {\gamma_t}{\eta_t}V_X(x,x_{t+1}) + \frac {\gamma_t}{\tau_t}V_Y(y, y_t) - \frac {\gamma_t}{\tau_t}V_Y(y, y_{t+1}) \\
				& +\gamma_t\langle K(x_{t+1}-x_{t}), y-y_{t+1}\rangle - \gamma_{t-1}\langle K(x_t-x_{t-1}), y-y_t\rangle\\
				& -\gamma_t\left(\frac{\alpha_X}{2\eta_t} - \frac{L_G}{2\beta_t}\right)\|x_{t+1}-x_t\|^2 - \frac {\alpha_Y\gamma_t}{2\tau_t}\|y_{t+1}-y_t\|^2  
				 -\gamma_{t-1}\langle K(x_t-x_{t-1}),y_{t}-y_{t+1}\rangle\\
				&\ - \gamma_t\langle\Delta_{x,G}^i+\Delta_{x,K}^i, x\tn-x\rangle - \gamma_t\langle\Delta_y^i, y\tn - y\rangle, \ \
				\forall z \in Z.
			\end{aligned}
		\end{equation}		
		By Cauchy-Schwartz inequality and \eqref{eqnIncStep}, for all $p\in (0, 1)$,
		\begin{align}
			\label{eqnCSS}
			\begin{aligned}
			&\ -\gamma_{t-1}\langle K(x_t-x_{t-1}),y_{t}-y_{t+1}\rangle
			\leq \gamma_{t-1}\|K(x_t-x_{t-1})\|_*\|y_t-y_{t+1}\|\\
			\leq &\ L_K\gamma_{t-1}\|x_t-x_{t-1}\|\|y_t-y_{t+1}\| \leq \frac{L_K^2\gamma_{t-1}^2\tau_{t}}{2p\alpha_Y\gamma_t}\|x_t-x_{t-1}\|^2+\frac{p\alpha_Y\gamma_t}{2\tau_{t}}\|y_{t}-y_{t+1}\|^2
			\\
			\leq & \frac{L_K^2\gamma_{t-1}\tau_{t-1}}{2p\alpha_Y}\|x_t-x_{t-1}\|^2+\frac{p\alpha_Y\gamma_t}{2\tau_{t}}\|y_{t}-y_{t+1}\|^2.
			\end{aligned}
		\end{align}	
		
		By \eqref{eqnIneqBetaGamma}, \eqref{eqnLambda}, \eqref{eqnACPDSRecursion} and \eqref{eqnCSS}, we can develop the following recursion for Algorithm \ref{algACPDS}:
		\begin{equation*}
			\begin{aligned}
				&\ (\beta_{t+1}-1)\gamma_{t+1}Q(z\ag[1], z) - (\beta_t-1)\gamma_tQ(z\ag , z)		
				= \ \beta_t\gamma_tQ(z\ag[1], z) - (\beta_t-1)\gamma_tQ(z\ag , z)\\
				\leq &\ \frac {\gamma_t}{\eta_t}V_X( x, x_t) - \frac {\gamma_t}{\eta_t}V_X( x,x_{t+1}) + \frac {\gamma_t}{\tau_t}V_Y( y, y_t) - \frac {\gamma_t}{\tau_t}V_Y( y, y_{t+1}) \\
				& +\gamma_t\langle K(x_{t+1}-x_{t}),  y-y_{t+1}\rangle - \gamma_{t-1}\langle K(x_t-x_{t-1}),  y-y_t\rangle\\
				& -\gamma_t\left(\frac{q\alpha_X}{2\eta_t} -  \frac{L_G}{2\beta_t}\right)\|x_{t+1}-x_t\|^2 + \frac{L_K^2\gamma_{t-1}\tau_{t-1}}{2p\alpha_Y}\|x_t-x_{t-1}\|^2 + \Lambda_t(x), \ \
				\forall z \in Z.
			\end{aligned}
		\end{equation*}	
	Applying the above inequality inductively and assuming that $x_0=x_1$, we obtain
			\begin{align}
				\label{eqnIneqSumS}
				\begin{aligned}
				&\ (\beta_{t+1}-1)\gamma_{t+1} Q(z^{ag}_{t+1}, z) - (\beta_1 - 1)\gamma_1Q(z^{ag}_1, z)
				\\		
				\leq & \, {\cal B}_t(z, z_{[t]})
 +\gamma_{t}\langle K(x_{t+1}-x_{t}),  y-y_{t+1}\rangle -\gamma_{t}\left(\frac{q\alpha_X}{2\eta_{t}}-\frac{L_G}{2\beta_{t}} \right)\|x_{t+1}-x_{t}\|^2
				\\
				&- \sum_{i=1}^{t-1}\gamma_i\left(\frac{q\alpha_X}{2\eta_i}  - \frac{L_G}{2\beta_i} - \frac{L_K^2\tau_i}{2p\alpha_Y} \right)\|x_{i+1} - x_i\|^2+\sum_{i=1}^{t}\Lambda_i(x), \ \
				\forall z \in Z.
				\end{aligned}
			\end{align}	
Relation \eqref{eqnQBoundGeneralS} then follows immediately from the above inequality,	\eqref{eqnIneqBetaGamma} and \eqref{eqnIneqLGLKS}.
	\end{proof}
\end{lemma}

\vgap

We also need the following technical result whose proof is based on Lemma 2.1 of \cite{nemirovski2009robust}.

\begin{lemma}
	\label{lemTech}
	Let $\eta_i, \tau_i$ and $\gamma_i$, $i = 1, 2, \ldots$, be given positive constants.
	For any $z_1\in Z$, if we define $z_1^v = z_1$ and 
	\begin{align}
		\label{eqnziv}
	z_{i+1}^v = \argmin{z=(x,y)\in Z}\left\{-\eta_i\langle \Delta^i_x, x\rangle -  \tau_i\langle\Delta^i_y, y\rangle + V(z, z_i^v) \right\},
	\end{align}
	then	
	\begin{align}
		\label{eqnTech}
		\sum_{i=1}^t\gamma_i\langle -\Delta_i, z_i^v - z\rangle 
		\leq 
		{\cal B}_t(z,z^v_{[t]}) +  \sum_{i=1}^t\frac{\eta_i\gamma_i}{2\alpha_X}\|\Delta^i_x\|_*^2 + \sum_{i=1}^t\frac{\tau_i\gamma_i}{2\alpha_Y}\|\Delta^i_y\|_*^2,
	\end{align}
	where $z^v_{[t]} := \{z^v_i \}_{i=1}^t$ and ${\cal B}_t(z,z^v_{[t]})$ is defined in \eqnok{def_CB}.

	\begin{proof}
	Noting that \eqref{eqnziv} implies $z_{i+1}^v = (x_{i+1}^v, y_{i+1}^v)$ where $x_{i+1}^v = \argmin{x=\in X}\left\{-\eta_i\langle \Delta^i_x, x\rangle  + V_X(x, x_i^v) \right\}$ and $y_{i+1}^v = \argmin{y\in Y}\left\{-  \tau_i\langle\Delta^i_y, y\rangle + V(y, y_i^v) \right\}$, 
	from Lemma 2.1 of \cite{nemirovski2009robust} we have
	\begin{align*}
		V_X(x, x^v_{i+1}) \leq & V_X(x, x_i^v) - \eta_i\langle \Delta^i_x, x - x_i \rangle + \frac{\eta_i^2\|\Delta^i_x\|^2_*}{2\alpha_X}, \\
		V_Y(y, y^v_{i+1}) \leq & V_Y(y, y_i^v) - \tau_i\langle \Delta^i_y, y - y_i \rangle + \frac{\tau_i^2\|\Delta^i_y\|^2_*}{2\alpha_Y}, 
	\end{align*}
	for all $i\geq 1$. Thus
	\begin{align*}
		\frac{\gamma_i}{\eta_i}V_X(x, x^v_{i+1}) \leq & \frac{\gamma_i}{\eta_i}V_X(x, x_i^v) - \gamma_i\langle \Delta^i_x, x - x_i \rangle + \frac{\gamma_i\eta_i\|\Delta^i_x\|^2_*}{2\alpha_X}, \\
		\frac{\gamma_i}{\eta_i}V_Y(y, y^v_{i+1}) \leq & \frac{\gamma_i}{\eta_i}V_Y(y, y_i^v) - \gamma_i\langle \Delta^i_y, y - y_i \rangle + \frac{\gamma_i\tau_i\|\Delta^i_y\|^2_*}{2\alpha_Y}.
	\end{align*}
	Adding the above two inequalities together, and summing up them from $i=1$ to $t$ we get
	\begin{align*}
		0 \leq {\cal B}_t(z, z^v_{[t]}) - \gamma_i\langle\Delta^i, z- z_i\rangle + \frac{\gamma_i\eta_i\|\Delta^i_x\|^2_*}{2\alpha_X} + \frac{\gamma_i\tau_i\|\Delta^i_y\|^2_*}{2\alpha_Y},
	\end{align*}
	so \eqref{eqnTech} holds. 
	\end{proof}
\end{lemma}

%
	
	\vgap
	
We are now ready to prove Theorem \ref{thmConvergenceBDS}.
	
\noindent	{\bf Proof of Theorem \ref{thmConvergenceBDS}}
Firstly, applying the bounds in \eqnok{bound_CB} and \eqnok{bound_inner}
to \eqref{eqnQBoundGeneralS}, we get
\begin{align}
		\label{eqnQBoundS}
		\begin{aligned}
			 \beta_t\gamma_t Q(z^{ag}_{t+1}, z) 					
			&\ \leq \frac {\gamma_t}{\eta_t}\Omega_X^2 - \frac{\gamma_{t}}{\eta_{t}}V_X( x, x_{t+1})
			+\frac{\gamma_t}{\tau_t}\Omega_Y^2 - \frac {\gamma_{t}}{\tau_{t}} V_Y( y, y_{t+1}) +\frac{\alpha_Y\gamma_{t}}{2\tau_{t}}\|y- y_{t+1}\|^2			\\
			& \ \ -\gamma_{t}\left(\frac{q\alpha_X}{2\eta_{t}}-\frac{L_G}{2\beta_{t}} - \frac{L_K^2\tau_t}{2\alpha_Y}\right)\|x_{t+1}-x_{t}\|^2 + \sum_{i=1}^{t}\Lambda_i(z) 
			\\
		 &\	\le  \frac{\gamma_t}{\eta_t}\Omega_X^2 + \frac{\gamma_t}{\tau_t}\Omega_Y^2 +  \sum_{i=1}^{t}\Lambda_i(z), \ \ 
		 \forall \, z \in Z.
		\end{aligned}
	\end{align}
By \eqref{eqnLambda}, we have
		\begin{align}
		\label{eqnDeltaBound}
		\begin{aligned}
			&\ \Lambda_i(z) = 	-\frac{(1-q)\alpha_X\gamma_i}{2\eta_i}\|x_{i+1}-x_i\|^2
			- \frac {(1-p)\alpha_Y\gamma_i}{2\tau_i}\|y_{i+1}-y_i\|^2 +\gamma_i\langle\Delta^i,  z-z_{i+1}\rangle\\
			= &\ -\frac{(1-q)\alpha_X\gamma_i}{2\eta_i}\|x_{i+1}-x_i\|^2
			- \frac {(1-p)\alpha_Y\gamma_i}{2\tau_i}\|y_{i+1}-y_i\|^2 
			+ \gamma_i\langle\Delta^i,  z_i-z_{i+1}\rangle +\gamma_i\langle\Delta^i,  z-z_{i}\rangle\\
			\le &\ \frac{\eta_i\gamma_i}{2(1-q)\alpha_X}\|\Delta_x^i\|_*^2 + \frac{\tau_i\gamma_i}{2(1-p)\alpha_Y}\|\Delta_y^i\|_*^2\
			 +\gamma_i\langle \Delta^i, z - z_{i}\rangle, 
		\end{aligned}
		\end{align}
where the last relation follows from Young's inequality. For all $i\geq 1$, letting $z^v_1 = z_1$, and $z^v_{i+1}$ as in \eqref{eqnziv}, we conclude from \eqnok{eqnDeltaBound} and Lemma~\ref{lemTech} that, $\forall z \in Z$,
\begin{align}
\label{def_Ut}
		\begin{aligned}
\sum_{i=1}^{t}\Lambda_i(z) \le &\ \sum_{i=1}^{t}
\left\{ \frac{\eta_i\gamma_i}{2(1-q)\alpha_X}\|\Delta_x^i\|_*^2 + \frac{\tau_i\gamma_i}{2(1-p)\alpha_Y}\|\Delta_y^i\|_*^2
+ \gamma_i\langle \Delta^i, z^v_i - z_{i}\rangle + \gamma_i\langle -\Delta^i, z^v_i - z\rangle\right\} 
\\
\le &\ \cB_t(z, z^v_{[t]}) + \underbrace{\frac 12\sum_{i=1}^{t}
\left\{ \frac{(2-q)\eta_i\gamma_i}{(1-q)\alpha_X}\|\Delta_x^i\|_*^2 + 
\frac{(2-p)\tau_i\gamma_i}{(1-p)\alpha_Y}\|\Delta_y^i\|_*^2
+ \gamma_i\langle \Delta^i, z^v_i - z_{i}\rangle \right\}}_{U_t},
	\end{aligned}
		\end{align}
where similar to \eqref{bound_CB} we have $\cB_t(z, z^v_{[t]})\leq \Omega_X^2\gamma_t/\eta_t + \Omega_Y^2\gamma_t/\tau_t$. 
Using the above inequality, \eqnok{def_gap}, \eqnok{boundness} and \eqnok{eqnQBoundS},
we obtain
\beq \label{eqnQTheoremS}
\beta_t\gamma_t g(z^{ag}_{t+1}) 
\le  \frac{2\gamma_t}{\eta_t} \Omega_X^2  +  \frac{2\gamma_t}{\tau_t} \Omega_Y^2 + U_t.
\eeq
Now it suffices to bound the above quantity $U_t$, both in 
expectation (part a)) and in probability (part b)).

We first show part a).  
		Note that by our assumptions on $\SO$, at iteration $i$ of Algorithm~~\ref{algACPDS},
		the random noises $\Delta^i$ are independent of $z_i$
		 and hence $\E[\langle\Delta^i, x-x_i\rangle]=0$.
		In addition, Assumption \ref{itmVar} implies that $\E[\|\Delta_{x}^i\|_*^2]\leq
		\sigma_{x,G}^2 + \sigma_{x,K}^2 = \sigma_{x}^2$
		(noting that $\Delta_{x,G}^i$ and $\Delta_{x,K}^i$ are indepdent at iteration $i$),
		and $\E[\|\Delta_y^i\|_*^2]\leq \sigma_y^2$. Therefore,
\begin{align}
	\label{eqnEUt}
\E[U_t] \le\frac{1}{2} \sum_{i=1}^{t}
\left\{ \frac{(2-q)\eta_i\gamma_i \sigma_x^2}{(1-q)\alpha_X} + 
\frac{(2-p)\tau_i\gamma_i \sigma_y^2}{(1-p)\alpha_Y} \right\}.
\end{align}
Taking expectation on both sides of \eqnok{eqnQTheoremS} and using
the above inequality, we obtain \eqnok{eqnEQBoundBDS}.

We now show that part b) holds. Note that by our assumptions on $\SO$
and the definition of $z^v_i$, the sequences $\{\langle \Delta_{x,G}^i, x_i^v - x_i\rangle\}_{i \ge 1}$
is a martingale-difference sequence. By the well-known large-deviation
theorem for matrigale-difference sequence (e.g., Lemma 2 of \cite{lan2012validation}),
and the fact that
	\begin{align*}
				& \E[\exp\left\{\alpha_Y\gamma_i^2\langle\Delta_{x,G}^i,  x^v_i - x_i\rangle^2 / \left(2\gamma_i^2\Omega_Y^2\sigma_{x, G}^2 \right) \right\} ]
				\leq \E[\exp\left\{\alpha_Y\|\Delta_{x,G}^i\|_*^2\| x^v_i - x_i\|^2 / \left(2\Omega_Y^2\sigma_{x, G}^2 \right) \right\} ]
				\\
				\leq & \E[\exp\left\{\|\Delta_{x,G}^i\|_*^2V( x^v_i,  x_i) / \left(\Omega_Y^2\sigma_{x, G}^2 \right) \right\} ]
				\leq  \E[\exp\left\{\|\Delta_{x,G}^i\|_*^2/ \sigma_{x, G}^2 \right\} ] \leq \exp\{1\},
			\end{align*}
we conclude that
\[
\begin{array}{l}
\prob\left\{\sum_{i=1}^{t}\gamma_i\langle\Delta_{x,G}^i,x^v_i - x_i\rangle>\lambda\cdot \sigma_{x, G} \Omega_X\sqrt{\frac 2{\alpha_X}\sum_{i=1}^{t}\gamma_i^2 }  \right\}\leq \exp\{-\lambda^2/3 \}, \forall \lambda >0.
\end{array}
\]
By using a similar argument, we can show that, $\forall \lambda > 0$,
\[
\begin{array}{l}
\prob\left\{\sum_{i=1}^{t}\gamma_i\langle\Delta_y^i, y^v_i - y_i\rangle>\lambda\cdot \sigma_y\Omega_Y\sqrt{\frac 2{\alpha_Y}\sum_{i=1}^{t}\gamma_i^2 }  \right\}\leq \exp\{-\lambda^2/3 \},\\
\prob\left\{\sum_{i=1}^{t}\gamma_i\langle\Delta_{x,K}^i,  x - x_i\rangle>\lambda\cdot \sigma_{x, K}\Omega_X\sqrt{\frac 2{\alpha_X}\sum_{i=1}^{t}\gamma_i^2 }  \right\}\leq \exp\{-\lambda^2/3 \}.
\end{array}
\]
Using the previous three inequalities and the fact that $\sigma_{x,G}+\sigma_{x,K} \le  \sqrt{2 \sigma_x}$,
we have, $\forall \lambda > 0$,
\beq \label{eqnProb1}
\begin{array}{l}
\prob\left\{\sum_{i=1}^{t}\gamma_i\langle\Delta^i, z^v_i - z_i\rangle>\lambda
\left[\frac{\sqrt{2} \sigma_{x} \Omega_X}{\sqrt{\alpha_X}} + \frac{\sigma_y \Omega_Y}{\sqrt{\alpha_Y}} \right]
 \sqrt{2\sum_{i=1}^{t}\gamma_i^2 }  \right\} \le\\
\prob\left\{\sum_{i=1}^{t}\gamma_i\langle\Delta^i, z^v_i - z_i\rangle>\lambda
\left[\frac{(\sigma_{x,G}+\sigma_{x,K}) \Omega_X}{\sqrt{\alpha_X}} + \frac{\sigma_y \Omega_Y}{\sqrt{\alpha_Y}} \right]
 \sqrt{2\sum_{i=1}^{t}\gamma_i^2 }  \right\}\leq 3\exp\{-\lambda^2/3 \}.
\end{array}
\eeq
	
Now let $S_i := (2-q)\eta_i\gamma_i/[(1-q)\alpha_X]$
and $S:=\sum_{i=1}^{t}S_i$. By the convexity of exponential function, we have
	\[
			\begin{array}{l}
				\E\left[\exp\left\{\frac 1S\sum_{i=1}^{t}S_i{\|\Delta_{x,G}^i\|_*^2}/{\sigma_{x,G}^2} \right\}\right] \leq \E\left[\frac 1S\sum_{i=1}^{t}S_i\exp\left\{\|\Delta_{x,G}^i\|_*^2/\sigma_{x,G}^2 \right\}\right]\leq\exp\{1\}.
			\end{array}
			\]
where the last inequality follows from Assumption \ref{itmLT}. Therefore, by Markov's inequality, for all $\lambda>0$,
			\[
			\begin{array}{l}
				\prob\left\{ \sum_{i=1}^{t} \frac{(2-q)\eta_i\gamma_i}{(1-q)\alpha_X}\|\Delta_{x,G}^i\|_*^2 > 
				(1+\lambda)\sigma_{x,G}^2 \sum_{i=1}^{t} \frac{(2-q)\eta_i\gamma_i}{(1-q)\alpha_X}\right\} \\
				=  \prob\left\{\exp\left\{\frac 1S\sum_{i=1}^{t}S_i{\|\Delta_y^i\|^2}/{\sigma_y^2} \right\}\geq \exp\{1+\lambda\} \right\}
				\leq 				\exp\{-\lambda\}.
			\end{array}
\]
Using an similar argument, we can show that
\[
\begin{array}{l}
\prob\left\{ \sum_{i=1}^{t}\frac{(2-q)\eta_i\gamma_i}{(1-q)\alpha_X}\|\Delta_{x,K}^i\|_*^2 >
(1+\lambda)\sigma_{x,K}^2 \sum_{i=1}^{t} \frac{(2-q)\eta_i\gamma_i}{(1-q)\alpha_X}\right\} \leq 	\exp\{-\lambda\},\\
\prob\left\{ \sum_{i=1}^{t} \frac{(2-p)\tau_i\gamma_i}{(1-p)\alpha_Y}\|\Delta_{y}^i\|_*^2 >
(1+\lambda)\sigma_{y}^2 \sum_{i=1}^{t} \frac{(2-p)\tau_i\gamma_i}{(1-p)\alpha_Y}\right\} \leq 	\exp\{-\lambda\}.
\end{array}
\]
Combining the previous three inequalities, we obtain
\beq \label{eqnProb2}
\begin{array}{l}
\prob\left\{ \sum_{i=1}^{t}\frac{(2-q)\eta_i\gamma_i}{(1-q)\alpha_X}\|\Delta_{x}^i\|_*^2
+  \sum_{i=1}^{t}\frac{(2-p)\tau_i\gamma_i}{(1-p)\alpha_Y}\|\Delta_{y}^i\|_*^2 > \right.\\
\left. (1+\lambda)\left[\sigma_{x}^2 \sum_{i=1}^{t} \frac{(2-q)\eta_i\gamma_i}{(1-q)\alpha_X}
+ \sigma_{y}^2 \sum_{i=1}^{t} \frac{(2-p)\tau_i\gamma_i}{(1-p)\alpha_Y}\right]
\right\} \leq 3	\exp\{-\lambda\},
\end{array}
\eeq
Our result now follows directly from \eqnok{def_Ut}, \eqnok{eqnQTheoremS}, \eqnok{eqnProb1} and \eqnok{eqnProb2}.
\endproof

\vgap

In the remaining part of this subsection, our goal is to prove Theorem~\ref{thmConvergenceUBS}, which describes the convergence rate of Algorithm \ref{algACPDS} when $X$ and $Y$ are both unbounded. Similar as proving Theorem \ref{thmConvergenceUB}, first we specialize the result of Lemma \ref{lemQBoundGeneralS} under \eqref{eqnIneqBetaGamma}, \eqref{eqnFlatStep} and \eqref{eqnIneqLGLKS}. The following lemma is analogous to Lemma \ref{lemEucl}.

	\begin{lemma}
		\label{lemEuclS}
		Let $\hat z = (\hat x, \hat y)\in Z$ be a saddle point of \eqref{eqnDSPP}. If $V_X(x, x_t)=\|x - x_t\|^2/2$ 
		and $V_Y(y, y_t)=\|y - y_t\|^2/2$ in Algorithm \ref{algACPDS}, and the parameters $\beta_t, \theta_t, \eta_t$
		and $\tau_t$ satisfy  \eqref{eqnIneqBetaGamma}, \eqref{eqnFlatStep} and \eqref{eqnIneqLGLKS}, then
		\begin{enumerate}
			\renewcommand{\theenumi}{(\alph{enumi})}
			\itemeqn
			\begin{align}
				\label{eqnIneqSadPtS}
				\begin{aligned}
				&\ \|\hat x - x\tn\|^2 + \|\hat x - x\tn^v\|^2 + \frac{\eta_t(1-p)}{\tau_t}\|\hat y - y\tn\|^2 + \frac{\eta_t}{\tau_t}\|\hat y - y\tn^v\|^2 
				\\
				\leq &\ 2\|\hat x - x_1\|^2 + \frac{2\eta_t}{\tau_t}\|\hat y - y_1\|^2 + \frac{2\eta_t}{\gamma_t}U_t, \text{ for all }t\geq 1,
				\end{aligned}
			\end{align}
			where $(x\tn^v$, $y\tn^v)$ and $U_t$ are defined in \eqref{eqnziv} and \eqref{def_Ut}, respectively.
			\item  
			\begin{equation}
				\label{eqnModifiedQS}
				\tilde g(z_{t+1}^{ag}, v\tn)\leq \frac 1{\beta_t\eta_t}\|x\ag[1] - x_1\|^2 + \frac 1{\beta_t\tau_t}\|y\ag[1] - y_1\|^2 + \frac{1}{\beta_t\gamma_t}U_t =: \delta_{t+1}, \text{ for all } t\geq 1,
			\end{equation}
			where $\tilde g(\cdot, \cdot)$ is defined in \eqref{def_gapt} and
			\begin{align}
				\label{eqnvS}
				v_{t+1} &=\left(\frac 1{\beta_t\eta_t}(2x_1 - x\tn - x\tn^v) , \frac 1{\beta_t\tau_t}(2y_1 - y\tn - y\tn^v) + \frac 1{\beta_t} K(x\tn - x_t) \right).
			\end{align}
		\end{enumerate}		
		\begin{proof}
			Apply \eqref{eqnIneqLGLKS}, \eqref{eqnLastProdp} and \eqref{def_Ut}  to \eqref{eqnQBoundGeneralS} in Lemma \ref{lemQBoundGeneralS}, we get
			\begin{align*}
				\beta_t\gamma_tQ(z\ag[1], z) \leq \bar \cB(z, z_t) + \frac{p\gamma_t}{2\tau_t}\|y - y_t\|^2 + \bar \cB(z, z^v_t) + U_t,
			\end{align*}			
			where $\bar {\cal B}(\cdot, \cdot)$ is defined as
			\begin{align*}
				\bar {\cal B}(z, \tilde z) := &\ \frac {\gamma_t}{2\eta_t}\|x - x_1\|^2  - \frac{\gamma_t}{2\eta_{t}}\| x- \tilde x\|^2 +\frac{\gamma_t}{2\tau_t}\|y- y_1\|^2  - \frac {\gamma_t}{2\tau_{t}}\|y- \tilde y\|^2. 
			\end{align*}
			thanks to \eqref{eqnFlatStep}. Now letting $z = \hat z$, and noting that $Q(z\ag[1], \hat z)\geq 0$, we get \eqref{eqnIneqSadPtS}.
			
		On the other hand, if we only apply \eqref{eqnIneqLGLKS} and \eqref{def_Ut} to \eqref{eqnQBoundGeneralS} in Lemma \ref{lemQBoundGeneralS}, then we get
		\begin{align*}	
		\begin{aligned}
			\beta_t\gamma_t Q(z^{ag}_{t+1}, z) 
			\leq &\ \bar \cB(z, z_t) +\gamma_t\langle K(x_{t+1}-x_{t}),  y-y_{t+1}\rangle + \bar \cB(z, z^v_t) + U_t.
		\end{aligned}
		\end{align*}				
		Apply \eqref{eqnFlatStep} and \eqref{eqnx1toxag} to $\cB(z, z_t)$ and $\bar \cB(z, z^v_t)$ in the above inequality, we get \eqref{eqnModifiedQS}.
	\end{proof}
	\end{lemma}

\vgap

	With the help of Lemma \ref{lemEuclS}, we are ready to prove Theorem \ref{thmConvergenceUBS}.
	
\noindent		{\bf Proof of Theorem \ref{thmConvergenceUBS}} 
		Let $\delta\tn$ and $v\tn$ be defined in \eqnok{eqnModifiedQS} and \eqnok{eqnvS}, respectively.
		Also let $C$ and $D$, respectively, be defined in \eqnok{eqnC} and \eqnok{eqnD}.
			It suffices to estimate $\E[\|v\tn\|]$ and $\E[\delta\tn]$. First it follows from
			\eqref{eqnFlatStep}, \eqref{eqnC} and \eqref{eqnEUt} that 
			\beq \label{bnd_U_t}
			\E[U_t] \le \frac{\gamma_t}{2 \eta_t} C^2.
			\eeq
			Using the above inequality, \eqnok{eqnFlatStep}, \eqnok{eqnD} and \eqnok{eqnIneqSadPtS}, we have
			\begin{align*}
				\E[\|\hat x - x\tn\|^2] \leq 2D^2 + C^2 \textrm{ and \ } \E[\|\hat y - y\tn\|^2] \leq (2D^2 + C^2){\frac{\tau_1}{\eta_1(1-p)}},
			\end{align*}
which, by Jensen's inequality, then imply that
			\begin{align*}
				\E[\|\hat x - x\tn\|] \leq \sqrt{2D^2 + C^2} \textrm{ and \ } \E[\|\hat y - y\tn\|^2] \leq \sqrt{2D^2 + C^2}\sqrt{\frac{\tau_1}{\eta_1(1-p)}}.
			\end{align*}
			Similarly, we can show that
			\begin{align*}
				\E[\|\hat x - x\tn^v\|] \leq \sqrt{2D^2 + C^2} \textrm{ and \ } \E[\|\hat y - y\tn^v\|^2] \leq \sqrt{2D^2 + C^2}\sqrt{\frac{\tau_1}{\eta_1}}.
			\end{align*}			
			Therefore, by \eqref{eqnvS} and the above four inequalities, we have
			\begin{align*}
				&\ \E[\|v_{t+1}\|] 
				\\
				\leq &\ \E\left[\frac 1{\beta_t\eta_t}\left(\|x_1 - x\tn\| + \|x_1 - x\tn^v\|\right) + \frac 1{\beta_t\tau_t}\left(\|y_1 - y\tn\|+\|y_1 - y\tn^v\|\right) + \frac {L_K}{\beta_t} \|x\tn - x_t\|\right]
				\\
				\leq &\ \E\left[\frac 1{\beta_t\eta_t}\left(2\|\hat x - x_1\| + \|\hat x - x\tn\| + \|\hat x - x\tn^v\|\right) \right.
				\\
				&\ \left.+ \frac 1{\beta_t\tau_t}\left(2\|\hat y - y_1\| + \|\hat y  - y\tn\| + \|\hat y  - y\tn^v\|\right) + \frac {L_K}{\beta_t} \left(\|\hat x - x\tn\| + \|\hat x - x_t\|\right)\right]
				\\
				\le &\ \frac {2\|\hat x - x_1\|}{\beta_t\eta_t} + \frac {2\|\hat y - y_1\|}{\beta_t\tau_t} + \sqrt{2D^2 + C^2} \left[\frac{2}{\beta_t\eta_t} + \frac{1}{\beta_t\tau_t}\sqrt{\frac{\tau_1}{\eta_1}}\left(\sqrt{\frac{1}{1-p}}+1\right) + \frac{2L_K}{\beta_t} \right],
			\end{align*}
			thus \eqref{eqnvboundS} holds.
	
			Now let us estimate a bound on $\delta\tn$. 
By \eqnok{eqnConvex_comb}, \eqref{eqnEUt}, \eqref{eqnIneqSadPtS} and \eqnok{bnd_U_t}, we have
			\begin{align*}
				&\ \E[\delta\tn]
				= \ \E\left[\frac 1{\beta_t\eta_t}\|x\ag[1] - x_1\|^2 + \frac 1{\beta_t\tau_t}\|y\ag[1] - y_1\|^2\right] + \frac{1}{\beta_t\gamma_t}\E[U_t]
				\\ 
				\leq &\ \E\left[\frac{2}{\beta_t\eta_t}\left(\|\hat x - x\ag[1]\|^2 + \|\hat x - x_1\|^2  \right) + \frac{2}{\beta_t\tau_t}\left(\|\hat y - y\ag[1]\|^2 + \|\hat y - y_1\|^2  \right)\right] + \frac{1}{2\beta_t\eta_t}C^2
				\\
				= &\ \E\left[\frac{1}{\beta_t\eta_t}\left(2D^2 + 2\|\hat x - x\ag[1]\|^2 + \frac{2\eta_t(1-p)}{\tau_t}\|\hat y - y\ag[1]\|^2 + \frac{2\eta_tp}{\tau_t}\|\hat y - y\ag[1]\|^2\right)\right] + \frac{1}{2\beta_t\eta_t}C^2
				\\
				\leq &\ \frac{1}{\beta_t\eta_t}\left[2D^2 + \frac{2}{\beta_t\gamma_t}\sum_{i=1}^{t}\gamma_i\left(\E\left[\|\hat x - x_{i+1}\|^2\right] + \frac{\eta_t(1-p)}{\tau_t}\E\left[\|\hat y - y_{i+1}\|^2\right]  + \frac{\eta_tp}{\tau_t}\E\left[\|\hat y - y_{i+1}\|^2\right]\right) + \frac{C^2}{2}\right]
				\\
				\leq &\ \frac{1}{\beta_t\eta_t}\left[2D^2 + \frac{2}{\beta_t\gamma_t}\sum_{i=1}^{t}\gamma_i\left(2D^2 + C^2 + \frac{\eta_tp}{\tau_t}\cdot\frac{\tau_1}{\eta_1(1-p)}(2D^2+C^2) \right)+ \frac{C^2}{2}\right] =
				\frac{1}{\beta_t\eta_t}\left(\frac{6-4p}{1-p}D^2 + \frac{5 - 3p}{2-2p}C^2 \right).
			\end{align*}
			Therefore \eqref{eqnErrUBDS} holds.
	
		\endproof

\section{Conclusion}
We present in this paper the APD method by incorporating a multi-step acceleration scheme into the primal-dual method in \cite{chambolle2011first}. 
We show that this algorithm can achieve the optimal rate of convergence for solving both deterministic and stochastic SPP. In particular,
the stochastic APD algorithm seems to be the first optimal algorithm for solving this important class of stochastic saddle-point problems
in the literature. For both deterministic and stochastic SPP, the developed APD algorithms can 
deal with either bounded or unbounded feasible sets as long as a saddle point of SPP exists.
In the unbounded case, the rate of convergence of the APD algorithms will depend on the distance from the initial point to the set of optimal solutions.

\bibliography{report}{}
\bibliographystyle{plain}

\end{document}